\documentclass[oneside,english]{amsart}
\usepackage[T1]{fontenc}
\usepackage[utf8]{inputenc}
\usepackage[a4paper]{geometry}
\usepackage{mathtools}
\usepackage{cases}
\usepackage{multirow}
\usepackage{booktabs}
\usepackage{color}
\usepackage{babel}
\usepackage{amstext}
\usepackage{amsthm}
\usepackage{amssymb}
\usepackage{colortbl}
\usepackage[binary-units]{siunitx}
\DeclareSIPrefix\micro{\text{\textmu}}{-3}
\usepackage{tikz}
\usepackage{pgfplotstable}
\usepackage{enumitem}
\usepackage{makecell}
\usepackage{bm}

\pgfplotsset{compat=1.16}

\pgfplotstableset{
	every head row/.style={before row=\toprule, after row=\midrule},
	every last row/.style={after row=\bottomrule},
}

\usepackage[unicode=true,pdfusetitle,
 bookmarks=true,bookmarksnumbered=false,bookmarksopen=true,
 breaklinks=false,pdfborder={0 0 1},backref=false,colorlinks=true]
 {hyperref}
\usepackage[noabbrev,capitalise,nameinlink]{cleveref}

\definecolor{badrow}{rgb}{1, 0.2, 0.2}
\definecolor{goodrow}{rgb}{0.4, 1, 0.4}

\makeatletter

\numberwithin{equation}{section}
\numberwithin{figure}{section}
\theoremstyle{plain}
\newtheorem{thm}{\protect\theoremname}
\theoremstyle{plain}
\newtheorem{prop}[thm]{\protect\propositionname}

\newcommand\restr[2]{{
  \left.\kern-\nulldelimiterspace 
  #1 
  \vphantom{\big|} 
  \right|_{#2} 
  }}

\makeatother

\newtheorem*{remark}{Remark}
\providecommand{\propositionname}{Proposition}
\providecommand{\theoremname}{Theorem}

\begin{document}
\title{Parallel kinetic schemes for conservation laws, with large time steps}
\author{Pierre Gerhard, Philippe Helluy, Victor Michel-Dansac, Bruno Weber}
\begin{abstract}
We propose a new parallel Discontinuous Galerkin method for the approximation of
hyperbolic systems of conservation laws. The method remains stable
with large time steps, while keeping the complexity of an explicit
scheme: it does not require the assembly and resolution of large linear
systems for the time iterations. The approach is based on a kinetic
representation of the system of conservation laws previously investigated
by the authors \cite{coulette2016palindromic,badwaik2018task,coulette2019high,drui2019analysis,gerhard2022unconditionally}.
In this paper, the approach is extended with a subdomain strategy
that improves the parallel scaling of the method on computers
with distributed memory.
\end{abstract}



\maketitle

\section{Introduction}
\label{sec:introduction}

The Discontinuous Galerkin (DG) method is generally used to approximate
hyperbolic systems of conservation laws, see for instance
\cite{bourdel1992resolution,fezoui2005convergence,hesthaven2007nodal,cockburn2012discontinuous}
and included references.
The DG method is well suited for parallel computations,
and it is often used in the context of large scale simulations.
However, the time step $\Delta t$ of the DG method
is limited by the Courant-Friedrichs-Lewy (CFL) condition,
which takes the form
\begin{equation*}
    \Delta t \leq K \frac h c,
\end{equation*}
where $h$ is the diameter of the smallest cell in the mesh,
$c$ is the maximal wave speed of the system of conservation laws,
and $K$ is a mesh-independent constant.
Unfortunately, this constant is often small, especially for high orders of approximation,
which impedes the efficiency of the method.
Such a restriction is often a problem because, in many applications,
the time step imposed by the CFL condition is much smaller than
time steps sufficient to ensure a good accuracy of the time integration.

To avoid such a restrictive condition,
one possibility is to construct time-implicit schemes,
which are free from a CFL condition,
but which involve having to solve a linear system.
However, the cost of inverting such a linear system can become prohibitive,
see for instance \cite{Catella2010,hochbruck2015implicit} and references therein.
One can also build locally implicit schemes,
but the associated CFL condition is still constrained
by the size of the smallest cell in the interface between the explicit and the implicit regions,
see \cite{bourdel1992resolution,dolean2010locally}.
Indeed, it is not easy to perform an automatic partitioning
whose interfaces would not contain small cells.
Yet another approach is to use a local time stepping strategy,
where more time steps are performed on the smallest cells than on larger cells,
see for instance \cite{dumbser2007arbitrary,muller2007fully,altmann2009local,diaz2009energy}.
In practice, the efficiency gained by using this method can be disheartening,
since there are usually many small cells,
on which lots of computations are still needed.
A common drawback of the above methods is that, in the context of meshes with cells of uniform size,
standard explicit methods will not be outperformed.

In the context of the finite difference method,
recent work \cite{ecer2000digital,gaffar2014explicit,Gaffar2015,yan2016explicit}
has proposed an unconditionally stable time integration method.
This method is based on the filtering of high frequencies in the solution
in order to recover a less restrictive CFL condition.
This approach relies on the computation of a few eigenvectors of the spatial operator,
associated to low frequency eigenvalues.

In \cite{gerhard2022unconditionally}, we have proposed
an unconditionally stable method whose complexity
(both in terms of computation time and storage) is in
$\mathcal{O}(n)$, with $n$ the number of degrees of freedom in the spatial approximation.
As a consequence, this method has the same complexity as an explicit scheme.
It relies on a vectorial kinetic interpretation of the system of conservation laws
based on \cite{bouchut1999construction,aregba2000discrete}.
The whole algorithm then reduces to the resolution of independent transport equations
coupled through relaxation source terms.
This whole system can be approximated implicitly while retaining the complexity of an explicit scheme.
Indeed, the implicit resolution of the transport equation with a DG scheme can be
performed with a downwind visiting of the mesh in the direction of the transport velocity,
and the implicit relaxation source terms can be applied locally.
In addition, the algorithm is parallelizable,
but dependencies in the computation reduce the parallel efficiency:
the scaling is not optimal when considering too large a number of threads.

In this paper, we propose an improvement of the aforementioned kinetic approach,
which relaxes the constraints in the parallel algorithm in order to improve the scaling.
The main idea is to perform a decomposition of the computational domain.
The transport equations are then solved with an iterative algorithm on each subdomain,
while the relaxation source terms are treated like in the non-partitioned case.
Hence, the main change compared to the non-partitioned case
consists in having to solve the transport equation multiple times, instead of once, per time iteration.
In fact, thanks to the structure of the DG transport solver,
we can prove that this iterative algorithm converges
to the solution of the fully implicit solver in at most three iterations,
under the non-restrictive CFL condition
\begin{equation*}
    \Delta t \leq \tilde K \frac H c,
\end{equation*}
where this time $H$ is the diameter of the smallest subdomain,
and $\tilde K$ still is a mesh-independent constant.
Moreover, this subdomain decomposition allows
a great improvement in the parallel scaling of the whole method,
because it relaxes dependencies in the transport solver.

The goal of this paper is to present this strategy and to validate this approach,
both by comparing it to established solvers
and by performing a large simulation a real-life situation
(the interaction of waves emitted by an antenna with a human body).
The paper is organized as follows.
In \cref{sec:kinetic_approximation}, we recall the vectorial kinetic approximation
of conservation and balance laws.
Then, the DG solver is presented in \cref{sec:implicit-dg}.
\Cref{sec:thread_based} explains the thread-based parallel algorithm,
while \cref{sec:subdomain_based} is devoted to the subdomain decomposition algorithm.
Finally, numerical results are proposed in \cref{sec:numerics}.
In \cref{sec:conclusion}, a conclusion concludes the paper.

\section{Kinetic approximation of first order conservations laws\label{sec:Kin-approx}}
\label{sec:kinetic_approximation}


\subsection{Kinetic approximation}

In this paper, we are interested in the numerical approximation of
a system of~$m$ conservation laws in dimension~$d$, governed by:
\begin{equation}
	\label{eq:conslaw}
	\partial_{t}W+\sum_{i=1}^{d}\partial_{i}Q^{i}(W)=0,
\end{equation}
where the unknown is a vector $W(X,t)\in\mathbb{R}^{m}$ depending
on the space variable $X=(x_{1}, \ldots, x_{d}) \in \Omega \subset \mathbb{R}^d$
and on the time variable $t \in \mathbb{R}$.
For the partial derivatives, we have introduced the notation
\[
	\partial_{i}=\frac{\partial}{\partial x_{i}},\quad\partial_{t}=\frac{\partial}{\partial t}.
\]

We assume the system to be hyperbolic.
To introduce this property, let $N=(N_1,\ldots,N_d)\in\mathbb{R}^{d}$
be an arbitrary space direction. The flux in direction $N$ is then defined by
\[
	Q(W,N)=\sum_{i=1}^{d}N_{i}Q^{i}(W).
\]
The Jacobian matrix of the flux $d_{W}(Q(W,N))$ is then supposed to be diagonalizable
with real eigenvalues~$\lambda_{r}(W,N)$, for $r \in \{ 1, \ldots, m \}$.
By definition, this means that the system \eqref{eq:conslaw} is hyperbolic.

The numerical approximation of such systems is, in general, a difficult subject.
One of the difficulties is that explicit schemes are subject to restrictive
time step conditions. Implicit schemes do not suffer from time step
conditions but require solving large sets of linear equations. In
previous work (see \cite{gerhard2022unconditionally} and included
references), we have proposed a method, based on a kinetic approach,
that avoids these constraints. We now recall the principles of this
kinetic representation.

We consider a set of $d+1$ kinetic velocities $V_{k} \in \mathbb{R}^d$,
$k\in \{ 0, \ldots, d \}$, associated to \textbf{vectorial} kinetic functions
$F_{k}(W)\in\mathbb{R}^{m}$.
Additional kinetic velocities and function could be introduced,
but we here the minimum number of velocities, $d+1$, is considered for simplicity.
The macroscopic data and the kinetic data are related by
\[
	W=\sum_{k=0}^d F_{k}.
\]
We also define ``Maxwellian'' equilibrium functions $M_{k}(W)\in\mathbb{R}^{m}$
such that
\[
	W=\sum_{k=0}^d M_{k}(W).
\]
The kinetic BGK representation then is given by transport equations with
relaxation source terms \cite{bouchut1999construction,aregba2000discrete}:
\begin{equation}
	\label{eq:kinlaw}
	\forall k\in \{ 0, \ldots, d \},
	\quad
	\partial_{t}F_{k}+V_{k}\cdot \nabla_{X}F_{k}=\frac{1}{\tau}\big(M_{k}(W)-F_{k}\big),
\end{equation}
with $\tau$ a (small) relaxation time.

When the relaxation time $\tau$ goes to $0^{+}$, the kinetic model \eqref{eq:kinlaw}
is formally equivalent to the initial system of conservation laws
\eqref{eq:conslaw} as long as
\begin{equation}
	\label{eq:maxwell_def}
	W=\sum_{k=0}^d M_{k}(W)
	\text{\quad and \quad}
	\forall i \in \{ 1, \ldots, d \},
	\quad
	\sum_{k=0}^d V_{k}^{i}M_{k}(W)=Q^{i}(W).
\end{equation}
Equations \eqref{eq:maxwell_def} constitute a linear system of size
$m(d+1)\times m(d+1)$ whose unknowns are the Maxwellian functions $M_k(W)$.
One expects this system to have a unique solution
as soon as the set of kinetic velocities is well chosen.

Theoretical arguments show that the formal limit described above actually is the true limit,
provided a so-called sub-characteristic condition
is satisfied \cite{bouchut1999construction,aregba2000discrete}.
This condition states that the kinetic velocities have to be greater
than the largest wave speed of the underlying hyperbolic system:
\[
	\forall k \in \{ 0, \ldots, d \}, \,
	\forall W \in \mathbb{R}^m, \,
	\forall N \in \mathbb{R}^d,
	\quad
	\left\|V_{k}\right\| > \max_{r}\left|\lambda_{r}(W,N)\right|.
\]

We now present the algorithm we use in practice to solve the kinetic equations \eqref{eq:kinlaw}.

\subsection{Kinetic algorithm}
\label{sec:kinetic_algorithm}

In practice, directly solving the BGK system \eqref{eq:kinlaw} is difficult.
It is usually better to split the equations into a transport step and a relaxation step.

At the initial time, we start with initial data $W(\cdot,0)$.
We have to choose kinetic vectors $F_{k}(\cdot,0)$ such that $W=\sum_{k}F_{k}$.
Obviously, this choice is not unique. A natural choice is to take
$F_{k}(\cdot,0)=M_{k}(W(\cdot,0))$ for all $k \in \{ 0, \ldots, d \}$.

Now, at each time step, to go from time $t$ to time $t + \Delta t$,
we adopt the following kinetic algorithm.
\begin{enumerate}
	\item Start with given kinetic data $F_{k}(\cdot,t)$ for all $k \in \{ 0, \ldots, d \}$:
	thus $W(\cdot,t)=\sum_{k}F_{k}(\cdot,t)$.
	\item \label{enu:transport-step}
          For each $k$ in $\{ 0, \ldots, d \}$,
	      solve the free transport equation
          $\partial_{t}F_{k}+V_{k} \cdot  \nabla_{X}F_{k}=0$
	      for a duration of $\Delta t$. In the continuous description of the
	      scheme, the free transport step can be solved exactly. It is given
	      by shift operations
	      \begin{equation}
		      \label{eq:shift}
		      \forall k \in \{ 0, \ldots, d \},
		      \quad
		      F_{k}(X,t+\Delta t^{-})=F_{k}(X-\Delta t V_{k},t).
	      \end{equation}
	      In practice we prefer to approximate this step by a Discontinuous
	      Galerkin solver.
          It is described below, in \Cref{sec:implicit-dg}.
	\item \label{enu:source-step}
          Define
	      \[
		      W(\cdot,t+\Delta t^{-})=\sum_{k=0}^d F_k(\cdot,t+\Delta t^{-}),
	      \]
	      and take
	      \[
		      W(\cdot,t + \Delta t^{+})=W(\cdot,t+\Delta t^{-}) \eqqcolon W(\cdot,t + \Delta t).
	      \]
	\item Apply a relaxation, with parameter $\omega \in [1, 2]$:
	      \begin{equation}
		      \label{eq:over_relaxation}
		      \forall k \in \{ 0, \ldots, d \},
		      \quad
		      F_{k}(\cdot,t + \Delta t^{+})=\omega M_{k}(W(\cdot,t + \Delta t^+))+(1-\omega)F_{k}(\cdot,t+\Delta t^{-}).
	      \end{equation}
\end{enumerate}

\begin{remark}
	Note that, even though $F_{k}(\cdot,0)=M_{k}(W(\cdot,0))$ at the initial time $t = 0$,
	this is no longer the case for subsequent iterations.
	Indeed, this would not be desirable, as it would lead to a scheme
	with first order accuracy in time.
\end{remark}

Let us point out that the approximation of the conservative data is
continuous in time:
\[
	W(\cdot,n\Delta t^{-})=W(\cdot,n\Delta t^{+}),
\]
while the approximation of the kinetic data is discontinuous at times
$t_{n}=n\Delta t$: in general,
\[
	F_{k}(\cdot,n\Delta t^{-})\neq F_{k}(\cdot,n\Delta t^{+}).
\]
In the presence of source terms on the conservation law,
the conservative data are no longer continuous in time.
The treatment of source terms is discussed in \Cref{subsec:source-term-management},
the next section.

In the relaxation step \eqref{eq:over_relaxation},
the parameter $\omega$ plays an important role. A natural
choice would be to take $\omega=1$. This corresponds to a projection
of the kinetic data on the Maxwellian state at the end of each time
step. This choice presents many interesting features: it leads to
an entropy dissipative, first order scheme. In addition, it is unconditionally
stable with respect to the time step $\Delta t$. It enters the large
category of kinetic schemes. It has been observed a long time ago that
these schemes are free of CFL conditions,
see for instance~\cite{brenier1984averaged,perthame90}.
However, this interesting property is rarely exploited in practical applications.
Another choice corresponds to taking $\omega = 2$,
leading to an over-relaxation procedure.
This choice, and its consequences, are described in \Cref{sec:equivalent_equation}.

\subsection{Handling source terms\label{subsec:source-term-management}}

This whole method can be extended to balance laws,
i.e., conservation laws with a source term, which take the form
\begin{equation*}
	\label{eq:conslaw-source}
	\partial_{t}W+\sum_{i=1}^{d}\partial_{i}Q^{i}(W)=S(W).
\end{equation*}
We refer, for instance, to \cite{coulette2019high,gerhard2022unconditionally},
where the kinetic equations are given by
\begin{equation*}
	\label{eq:kinlaw_with_source}
	\forall k\in \{ 0, \ldots, d \},
	\quad
	\partial_{t}F_{k}+V_{k}\cdot \nabla_{X}F_{k} = G_k + \frac{1}{\tau}\left(M_{k}(W)-F_{k}\right),
\end{equation*}
with $G_k = (\nabla_W M_k(W)) S(W)$ for all $k\in \{ 0, \ldots, d \}$.

In presence of source terms, \Cref{enu:source-step} of the kinetic
algorithm is modified as follows.

\begin{enumerate}
	\item[(3)]
		Define
		\[
			W(\cdot,\Delta t^{-})=\sum_{k}F_{k}(\cdot,\Delta t^{-}).
		\]
		Solve the differential equation
		\begin{equation}
			\partial_{t}U(\cdot,t)=S(U(\cdot,t)),\label{eq:source-eqdiff}
		\end{equation}
		with the initial condition
		\[
			U(\cdot,0)=W(\cdot,\Delta t^{-}).
		\]
		Finally, take
		\[
			W(\cdot,\Delta t^{+})=U(\cdot,\Delta t).
		\]
\end{enumerate}

In practice, the differential equation \eqref{eq:source-eqdiff}
is discretized with a Crank-Nicolson scheme, which reads
\begin{equation}
    \label{eq:source_term_application_to_W}
	\frac{W(\cdot,\Delta t^{+})-W(\cdot,\Delta t^{-})}{\Delta t}=\frac{S(W(\cdot,\Delta t^{+}))+S(W(\cdot,\Delta t^{-}))}{2}.
\end{equation}
For more details, we refer to \cite{gerhard2022unconditionally}.

\subsection{Equivalent equation}
\label{sec:equivalent_equation}

Because the first-order scheme is generally not accurate enough, it
is often better to consider the over-relaxed choice $\omega=2$. In
this case, the scheme becomes second-order accurate.
We briefly sketch the proof of this property.
For more details, we refer to \cite{coulette2019high,gerhard2022unconditionally}.

During the computations, one expects that,
for all $k \in \{0,\ldots,d\}$, $F_k \simeq M_k(W)$,
and therefore that $\sum_{k}V_{k}^{i}F_{k}\simeq Q^{i}(W)$.
We thus introduce the \textbf{approximate flux} $Z$ and the
\textbf{flux error} $Y$, defined as follows:
\begin{equation}
	\label{eq:def_Z_and_Y}
	Z^{i} \coloneqq \sum_{k=0}^d V_{k}^{i}F_{k},
	\text{\quad and \quad}
	Y^{i} \coloneqq Z^{i}-Q^{i}(W) = \sum_{k=0}^d V_{k}^{i} \left( F_{k} - M_k(W) \right).
\end{equation}
The whole kinetic algorithm is a functional operator $\mathcal{M}(\Delta t)$
that maps $(W(\cdot,0),Y^{i}(\cdot,0))$ to $(W(\cdot,\Delta t),Y^{i}(\cdot,\Delta t^{+}))$.
The operator $\mathcal{M}$ is made of (linear) shift operations and
(non-linear) local relaxations.
In the $(W,Y^{i})$ variables, the relaxation operation \eqref{eq:over_relaxation} simply reads, arguing \eqref{eq:def_Z_and_Y}:
\begin{equation*}
	\left\{
	\begin{aligned}
		\omega = 1
		 & \implies
		F_{k}(\cdot,\Delta t^{+})=M_{k}(W(\cdot,\Delta t)),                                                    \\
		 & \implies
		Y^{i}(\cdot,\Delta t^{+})
		=
		\sum_{\smash k=0}^{\smash d} V_{k}^{i} \left( F_{k}(\cdot,\Delta t^{+}) - M_k(W(\cdot,\Delta t)) \right)
		=
		\sum_{\smash k=0}^{\smash d} V_{k}^{i} \left( M_k(W(\cdot,\Delta t)) - M_k(W(\cdot,\Delta t)) \right)  \\
		 & \implies
		Y^{i}(\cdot,\Delta t^{+})
		=
		0 ,                                                                                                    \\
		\omega = 2  \vphantom{\dfrac 1 2}
		 & \implies
		F_{k}(\cdot,\Delta t^{+})=2 M_{k}(W(\cdot,\Delta t)) - F_{k}(\cdot,\Delta t^{-}),                      \\
		 & \implies
		Y^{i}(\cdot,\Delta t^{+})
		=
		\sum_{\smash k=0}^{\smash d} V_{k}^{i} \left( F_{k}(\cdot,\Delta t^{+}) - M_k(W(\cdot,\Delta t)) \right)
		=
		\sum_{\smash k=0}^{\smash d} V_{k}^{i} \left(M_k(W(\cdot,\Delta t) - F_{k}(\cdot,\Delta t^{-})\right), \\
		 & \implies Y^{i}(\cdot,\Delta t^{+})
		=
		- Y^{i}(\cdot,\Delta t^{-}).
	\end{aligned}
	\right.
\end{equation*}
Therefore, the choice $\omega = 2$ induces fast oscillations of the flux error.
For the forthcoming analysis, it is thus better to replace $\mathcal{M}$ with $\mathcal{M}\circ\mathcal{M}$.

In principle it is now easy, although tedious, to compute the equivalent
equation of the kinetic algorithm. It consists in computing a Taylor
expansion of
\[
	\frac{\mathcal{M}(\Delta t/2)-\mathcal{M}(-\Delta t/2)}{\Delta t},
\]
with respect to $\Delta t$, up to order $\mathcal{O}(\Delta t^{2})$.

During the calculation of the Taylor expansion, the term $X-\Delta t V_{k}$
in the shift operation \eqref{eq:shift} generates partial derivatives
in space. In addition, because of symmetries, when $\omega=2$, the
even-order terms of the expansion vanish. And finally, the relaxation
introduces non-linearities. We end up with a system of non-linear
partial differential equations of first order in $(W,Y^{i})$. The
calculations are tedious, but can be automated through a Computer
Algebra System.

We illustrate the results obtained for $d=1$, i.e., in one space dimension.
In this case, we have $d+1=2$ kinetic velocities.
We set $\lambda \in \mathbb{R}_+^*$ and we choose $V_{0}=-\lambda$ as well as $V_{1}=\lambda$.
In this one-dimensional case, $W=F_{0}+F_{1}$.
The equivalent equation for $\omega=2$ is then, in conservative variables and up to $\mathcal{O}(\Delta t^{2})$:
\[
	\partial_{t}W+\partial_{x}Q(W)=0.
\]
We indeed recover the desired conservation laws. It is also possible
to compute a second-order equivalent equation for $Y$:
\[
	\partial_{t}Y-d_{W}(Q(W))\partial_{x}Y=0.
\]

We observe that the system is hyperbolic and that the waves for $W$
and $Y$ move in opposite directions \cite{drui2019analysis}. We
emphasize that there is
no assumption of smallness of $Y$. In practice, we indeed observe
second-order accuracy, even when the initial $Y$ is of order $\mathcal{O}(1)$.
For $d=1$, the second-order expansion is not sufficient to analyze
the stability of the approximation. By analyzing the third-order term,
it is, however, possible to prove stability under a sub-characteristic
stability condition
\[
	\lambda\geq\max_{1\leq i\leq m}\left|\lambda_{i}(W)\right|,
\]
where $\lambda_{i}(W)$ are the eigenvalues of $d_{W}(Q(W))$, see \cite{drui2019analysis}.

We can also perform the calculations in the case $d=2$.
We then need $d+1=3$ kinetic velocities.
We can take
\[
	V_{k}=\begin{pmatrix}
		\cos \left( \dfrac{2k\pi}{3} \right) \\[\bigskipamount]
		\sin \left( \dfrac{2k\pi}{3} \right)
	\end{pmatrix}
\]
The equivalent equation on $W$ for $\omega=2$ at order $\mathcal{O}(\Delta t^{2})$
is, of course, the system of conservation laws \eqref{eq:conslaw}
\[
	\partial_{t}W+\partial_{1}Q^{1}(W)+\partial_{2}Q^{2}(W)=0.
\]
Setting $A^{i}(W)=d_{W}(Q^{i}(W))$, the equation for $Y$ is
\begin{equation*}
	\partial_{t}\begin{pmatrix}
		Y_{1} \vphantom{\dfrac 1 2} \\[\bigskipamount]
		Y_{2} \vphantom{\dfrac 1 2}
	\end{pmatrix}
	+
	\begin{pmatrix}
		\dfrac{\lambda}{2}I-A^{1}(W) & 0                   \\[\bigskipamount]
		-A^{2}(W)                    & -\dfrac{\lambda}{2}
	\end{pmatrix}
	\partial_{1}\begin{pmatrix}
		Y_{1} \vphantom{\dfrac 1 2} \\[\bigskipamount]
		Y_{2} \vphantom{\dfrac 1 2}
	\end{pmatrix}
	+ \\
	\begin{pmatrix}
		0                   & -\dfrac{\lambda}{2}I-A^{1}(W) \\[\bigskipamount]
		-\dfrac{\lambda}{2} & -A^{2}(W)
	\end{pmatrix}
	\partial_{2}\begin{pmatrix}
		Y_{1} \vphantom{\dfrac 1 2} \\[\bigskipamount]
		Y_{2} \vphantom{\dfrac 1 2}
	\end{pmatrix}
	=\mathcal{O}(\Delta t^{2}).
\end{equation*}
We can prove that the equivalent system is hyperbolic if $\lambda$ is large enough.
In this case, the sub-characteristic condition arises
from the analysis of the first order terms.
To get a finer bound, there exist more sophisticated analyses, based on the entropy:
we refer for instance to \cite{bouchut1999construction,dubois2014simulation}.

\section{Unconditionally stable DG approximations}
\label{sec:implicit-dg}


In this section, we recall how to construct an unconditionally stable
Discontinuous Galerkin (DG) approximation of the initial system of
conservation laws \eqref{eq:conslaw}.
The reader is referred to \cite{gerhard2022unconditionally},
where this procedure is explained in detail.

The kinetic algorithm presented in \Cref{sec:Kin-approx} relies
on transport steps and relaxation steps. The relaxation step is generally
easy to implement at each interpolation point of the approximation.
In addition, it is embarrassingly parallel.

The implementation difficulty of the kinetic algorithm
lies in the transport step, given by \Cref{enu:transport-step}.
In this step, the shift operation \eqref{eq:shift} consists in solving $(d+1) \times m$ transport equations of the form
\begin{equation}
	\label{eq:transport}
	\partial_{t}f + V \cdot \nabla f = 0.
\end{equation}
Note that $V$ represents one of the $(d+1)$ kinetic velocities
and that $f$ represents one of the $m$ dimensions of the $(d+1)$ kinetic unknowns.

If the computational domain has a simple shape and if the solution
is computed on a structured Cartesian grid, it is natural to solve
this transport equation by the characteristic method \eqref{eq:shift}.
With well-chosen time step $\Delta t$ and kinetic velocities $V_{k}$,
this approach leads to the so-called Lattice Boltzmann method,
see for instance \cite{baty:hal-02965967}, and included references.

In a domain $\Omega$ with a complex geometry, or discretized with an unstructured grid,
the characteristic method is no longer a good choice because it
leads to difficulties such as instabilities or loss of the conservation property.
In addition, the treatment of boundary conditions is not natural in
this framework. Instead, in the unstructured case, we prefer to rely on an approximation of the shift operation \eqref{eq:shift},
based on discretizing \eqref{eq:transport} in the DG framework.
For a general presentation of the DG approach, we refer to the book of Hesthaven and Warburton~\cite{hesthaven2007nodal}.
The idea to solve a kinetic BGK model with a DG approximation of the
transport step was already proposed in \cite{shi2003discontinuous}, but with
an explicit scheme. The novelty of our approach is to adopt an implicit
DG approximation, instead of an explicit one, in order to get rid
of the CFL condition. It turns out that the implicit scheme is not more
complicated to solve than the explicit scheme, because the matrix
of the implicit step is triangular and can thus be solved in an explicit fashion.
For the sake of completeness, we briefly describe the implicit DG
approach to approximate solutions to \eqref{eq:transport}. More details can be found
in \cite{badwaik2018task,coulette2019high,gerhard2022unconditionally}.

\subsection{DG scheme}

We consider an unstructured mesh $\mathcal{M}$ of the computational domain $\Omega$
made of tetrahedral cells.
On each cell, we define $n_n$ basis functions
$(\psi_{j}^{L}(X))_{j \in \{1,\dots,n_n\}}$.
The transported function $f$ is then approximated
in cell $L$ by a linear expansion on basis functions
\[
	f(x,n\Delta t)\simeq f_{L}^{n}(X)=\sum_{j}^{n_n} f_{L,j}^{n}\psi_{j}^{L}(X),\quad X \in L.
\]
The unknowns of the scheme are the coefficients $f_{L,j}^{n}$ of the linear expansion.

\begin{figure}[!ht]
	\centering
	\includegraphics{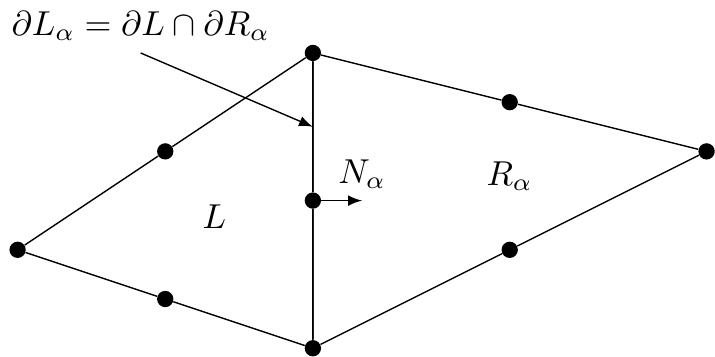}
	\caption{Notation for the Discontinuous Galerkin approximation.\label{fig:dg_geom}}
\end{figure}

We now write an implicit DG approximation scheme to compute the unknown coefficients $f_{L,j}^{n}$ at time $t^n$ from the known coefficients $f_{L,j}^{n-1}$ at time $t^{n-1}$.
For simplicity, we describe the case of an implicit first order Euler method.
The strategy can be extended to other more accurate schemes, such as the Crank-Nicolson scheme (which we use in practice) or DIRK (Diagonally Implicit Runge-Kutta) approaches,
see for instance \cite{alexander1977diagonally,KenCar2019}.
The DG scheme then reads as follows: for
each cell $L$ and each basis function $\psi_{i}^{L}$,
\begin{equation}
	\label{eq:dg_imp}
	\int_{L}\frac{f_{L}^{n}-f_{L}^{n-1}}{\Delta t}\psi_{i}^{L}
    -
    \int_{L} f_{L}^{n} \, V\cdot \nabla\psi_{i}^{L}
    +
    \sum_{\alpha = 1}^{n_f} \int_{\partial L_\alpha} \!\!%
    \Big(f_{L}^{n} \, (V\cdot N_\alpha)_{+} + f_{R_\alpha}^{n} \, (V\cdot N_\alpha)_{-} \Big)
    \psi_{i}^{L}=0.
\end{equation}
In this formula, $n_f$ denotes the number of faces of cell $L$
(for a tetrahedron, $n_f = 4$),
$\partial L_\alpha$ denotes the part of the boundary of $L$
where face $\alpha$ is located,
and $R_\alpha$ denotes the neighboring cell along $\partial L_\alpha$.
The situation is depicted in \Cref{fig:dg_geom}, in 2D for simplicity.
In addition, we use standard notation for the upwind numerical flux:
\[
	(V \cdot N_\alpha)_{+}=\max(V \cdot N_\alpha,0)
    \text{\quad and \quad}
    (V \cdot N_\alpha)_{-}=\min(V \cdot N_\alpha,0),
\]
where the vector $N_\alpha$ is the unit normal vector on
$\partial L_\alpha$ oriented from $L$ to $R_\alpha$.

\subsection{\label{subsec:downwind}Downwind algorithm}

The scheme is implicit and it seems that one would need to assemble
and solve a large linear system in order to compute $f_{L}^{n}$ from $f_{L}^{n-1}$.
However, we can exhibit an algorithm with explicit complexity, dubbed \emph{downwind algorithm},
that solves efficiently -- and in parallel -- the set of equations \eqref{eq:dg_imp}.
The method is described in detail
in~\cite{coulette2016palindromic,badwaik2018task,coulette2019high,gerhard2022unconditionally}.
In this manuscript, we only recall its major steps.

In the preprocessing phase, we construct a graph $G$ from the mesh $\mathcal{M}$.
Its nodes correspond to the cells of the mesh and its edges to the
faces between cells.
Each edge is then oriented with respect to the velocity $V $.
Between two nodes $L$ and $R$ (corresponding to two cells),
the edge is oriented from $L$ to $R$ if $V $ is oriented from $L$ to $R$,
i.e. if $V $ crosses the edge $\partial L \cap \partial R$ from $L$ to $R$
Because the velocity is constant, it is possible
to prove that the graph $G$ is direct and acyclic.
It can thus be sorted in topological order,
using Breadth-First Search
(see~\cite{gerhard2022unconditionally} for a comparison between
Breadth-First Search and Depth-First Search on such problems).
Note that this preprocessing phase is executed
only once at the beginning of the computations.

\begin{figure}[!ht]
		\centering
        \begin{tikzpicture}
            \node (left_graph) at (0,0) {
                \includegraphics{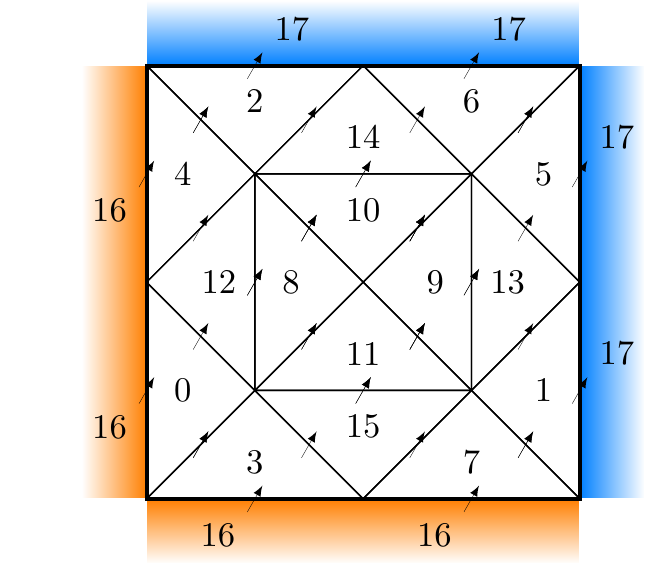}
            };
            \node [anchor=west, xshift=1cm] (right_graph) at (left_graph.east) {
                \includegraphics{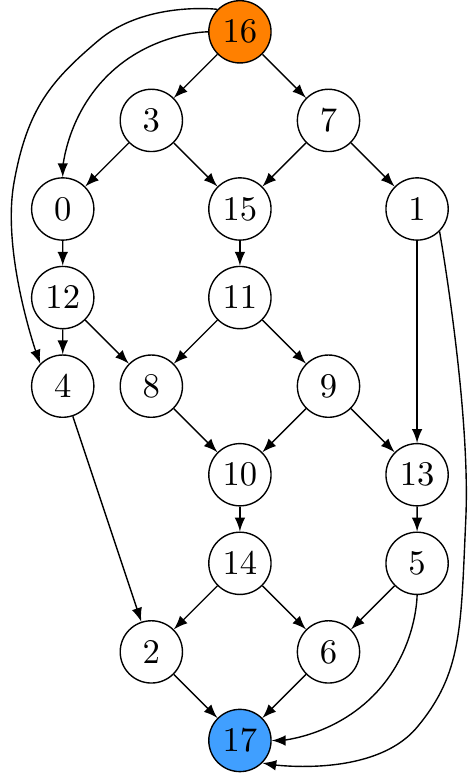}
            };
        \end{tikzpicture}
		\caption{%
            Example of a mesh $\mathcal{M}$ (left panel) and its associated graph $G$ (right panel).
            The nodes of the graph correspond to the cells of the mesh.
            Two additional, fictitious nodes are considered:
            the upwind node (in orange) and the downwind node (in blue).
            The solution can be explicitly computed by following
			a topological ordering of a Direct Acyclic Graph (DAG) using Breadth-First Search,
            e.g. 3, 7, 0, 15, 1, \textit{etc}.
            In addition, the parallel capabilities of the method are visible on the graph:
            first, cells $3$ and $7$ can be computed in parallel;
            then cells $0$, $15$ and $1$ can be computed in parallel, \textit{etc.}
        }
        \label{fig:dag}
	\end{figure}

In the main computation phase, corresponding to the time loop,
the linear system \eqref{eq:dg_imp} is then solved by
visiting the cells of the mesh in this topological order.
The algorithm is parallel and its storage can be optimized:
the solution can be replaced in memory during the computations,
see \Cref{fig:dag} for a example.
After solving the linear system, the relaxation step,
which is embarrassingly parallel, is applied.

\section{Thread-based numerical implementation}
\label{sec:thread_based}

We have implemented the kinetic algorithm in a code written in Rust.
Rust is a recent programming language oriented toward security and efficiency.
Most common bugs are avoided at compile time. For instance,
memory leaks, segmentation faults, uninitialized data and race conditions
are forbidden by the compiler. In addition, Rust proposes automatic
parallelization tools through the \texttt{rayon} library, based on
a work stealing strategy~\cite{matsakis2019rayon}.
This library is particularly well suited to the parallel implementation of
the downwind algorithm presented in \Cref{subsec:downwind}.
For more detail, we refer to \cite{gerhard2022unconditionally}.
Furthermore, meshes are generated using the \texttt{Gmsh} tool,
described in~\cite{geuzaine2009gmsh}.

The implementation presents a good scaling for a moderate number of
threads, as shown in \Cref{tab:multithread} from \cite{gerhard2022unconditionally}.

\begin{table}[!ht]
    \centering
    \caption{%
        Multithread efficiency of the downwind algorithm
        executed on a server with an Intel Xeon E5-2680 v3 processor
        ($24$ physical cores, \SI{2.50}{\giga\hertz}).
        The scalability is computed on coarse to fine meshes,
        with several refinement levels.
        We observe that the efficiency stalls at around~$\SI{60}{\percent}$.
    }
    \label{tab:multithread}
    \begin{tabular}{cccccccc}
        \toprule
        \multicolumn{2}{c}{refinement} & \multicolumn{2}{c}{it/s} & \multicolumn{2}{c}{\SI{}{\micro\second}/dof/it} &                                                       \\
        \cmidrule(lr){1-6}
        level                          & elements                 & serial                            & parallel & serial & parallel & scalability & heap     \\
        \cmidrule(lr){1-2} \cmidrule(lr){3-4} \cmidrule(lr){5-6} \cmidrule(lr){7-7} \cmidrule(lr){8-8}
        8                              & 1808                     & 72.58                             & 346.1    & 0.425  & 0.089    & 4.769       & \SI{11.85}{\mega\byte} \\
        16                             & 9199                     & 11.34                             & 102.2    & 0.569  & 0.063    & 9.012       & \SI{42.50} {\mega\byte} \\
        32                             & 56967                    & 1.698                             & 20.19    & 0.664  & 0.056    & 11.89       & \SI{266.5} {\mega\byte} \\
        48                             & 175138                   & 0.531                             & 7.753    & 0.718  & 0.049    & 14.60       & \SI{808.9} {\mega\byte} \\
        64                             & 386806                   & 0.236                             & 3.579    & 0.747  & 0.049    & 15.17       & \SI{1.777} {\giga\byte} \\
        72                             & 544030                   & 0.165                             & 2.531    & 0.765  & 0.050    & 15.34       & \SI{2.515} {\giga\byte} \\
        \bottomrule
    \end{tabular}
\end{table}

However, in \Cref{tab:multithread}, we also observe that the efficiency is not $\SI{100}{\percent}$.
Indeed, dependencies in the computations limit the parallel scaling of the downwind algorithm.
For instance, in the mesh from \Cref{fig:dag},
it is clear that launching more than three threads
is useless because the additional threads will have to wait
for computations to be finished before starting to work.
Similar behavior occurs for larger meshes,
whose parallel regions (the cells which can be treated in parallel) are,
on average, larger, but which contain small, efficiency-limiting parallel regions close to edges.

To address this issue, we propose in the next section a modification of the
downwind algorithm, both to improve the parallel efficiency of the method,
and to deal with the distributed-memory setting.

\section{Subdomain parallelism}
\label{sec:subdomain_based}


As explained above, the downwind algorithm is parallelized with a
work stealing thread-based algorithm. The parallel scaling is good
for a few threads,
but seems to be capped as the mesh becomes finer.

In order to provide better scaling capabilities, we now describe a subdomain strategy that
relaxes the computation dependencies. The main idea is to apply
the above time-implicit downwind algorithm in each subdomain, but
with a time-explicit coupling between the subdomains, so as to relax
the dependencies between regions. Because of the explicit coupling, it will become necessary
to apply an iterative algorithm to compute the approximate solution in
a stable way. The algorithm can be proved to converge in a finite
number of iterations. In most configurations, three iterations are
sufficient. Let us now describe the principles of this subdomain iterative
algorithm.

As in Section \ref{sec:implicit-dg}, the main task is the resolution
of the initial value problem for the transport equation:
\begin{equation*}
	\begin{dcases}
		\partial_{t}f+V \cdot \nabla f=0 & \text{for } X \in \Omega\times[0,\Delta t], \\
		f(X,0)=f^{0}(X)                & \text{for } X \in \Omega.
	\end{dcases}
\end{equation*}

\subsection{Iterative algorithm}

We assume that $\Omega$ is decomposed into a finite number of subdomains
$(\Omega_{i})_{i \in \{1, \ldots, n_{d}\}}$. To simplify the presentation,
we assume that $\Omega$ is either a periodic domain or the whole space domain,
in order to avoid having to describe the boundary conditions. However,
the approach is also valid when~$\partial\Omega\neq\emptyset$.

We then denote
by $f_{i}$ the restriction of $f$ to subdomain $\Omega_{i}$,
by $N_{i}(X)$ the outward normal vector on $\partial\Omega_{i}$,
by $\mathcal{N}(\Omega_i)$ the subdomains neighboring $\Omega_i$,
and by $\partial\Omega_{i}^{-}$ the upwind part of the boundary of
$\Omega_{i}$:
\[
	\partial\Omega_{i}^{-}=\left\{ X\in\partial\Omega_{i} \mid N_{i}(X)\cdot V <0\right\} .
\]

We initialize the algorithm by setting $f_{i}^{0}(X,t)=f_{i}^{0}(X)$.
Thus, the initial iteration does not depend on time.

We then propose an iterative algorithm to compute the successive
time-dependent iterations $f_{i}^{p}$ in subdomain~$\Omega_{i}$,
for $p \geq 1$.
To compute $f_{i}^{p}$ from $f_{i}^{p-1}$,
we solve the following time-dependent boundary value problems:
\begin{subnumcases}{\label{eq:ibvp}}
	\label{eq:ibvp1}
	\partial_t f_i^p + V  \cdot \nabla f_i^p = 0 & $\text{for } X \in \Omega_i \text{\, and \,} t \in (0, \Delta t)$,                                                   \\
	\label{eq:ibvp2}
	f_i^p(X, 0) = f_i^0(X)               & $\text{for } X \in \Omega_i$,                                                                                        \\
	\label{eq:ibvp3}
	f_i^p(X, t) = f_j^{p-1}(X, t)        & $\text{for } X \in \partial \Omega_i^- \cap \partial \Omega_j \text{, } \forall \Omega_j \in \mathcal{N}(\Omega_i)$.
\end{subnumcases}
We can then prove the following result.
\begin{prop}
	let $\mathcal{L}$ be the maximal subdomain diameter.
	Under the condition
	\[
		\Delta t\leq\frac{\mathcal{L}}{\left| V \right|},
	\]
	the above algorithm \eqref{eq:ibvp}
	converges to the exact solution in at most three iterations: $f_{i}^{3}=f_{i}$.
\end{prop}

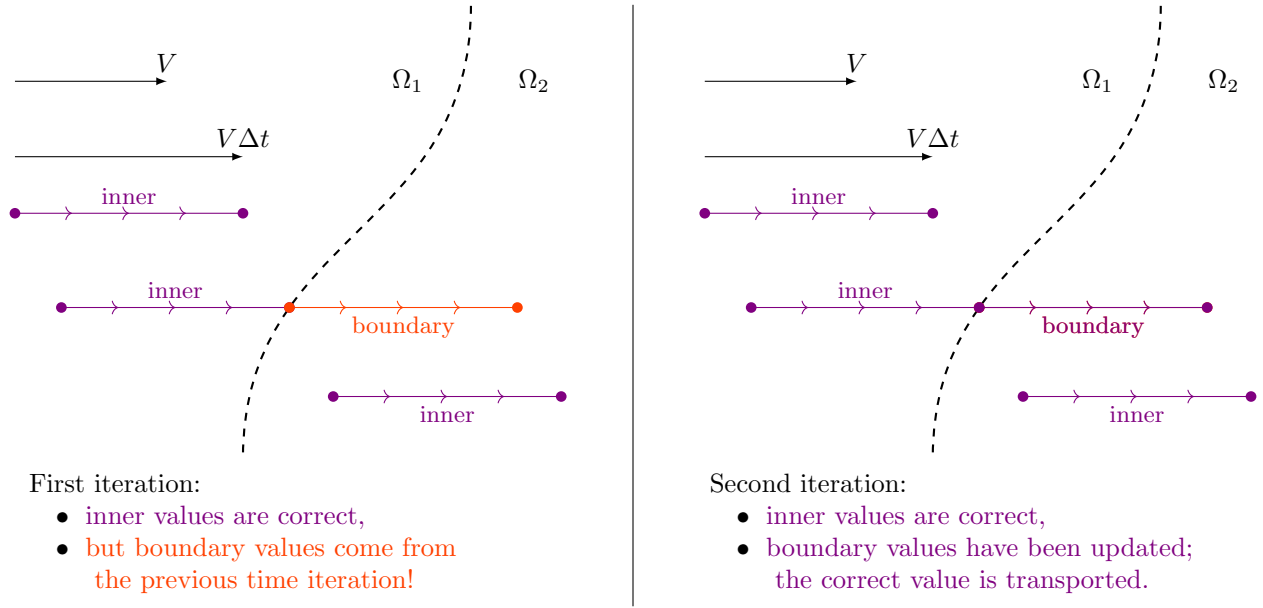
\begin{figure}[!ht]

	\centering

	\begin{tikzpicture}
		\node (left_graph) at (0, 0) {\includegraphics{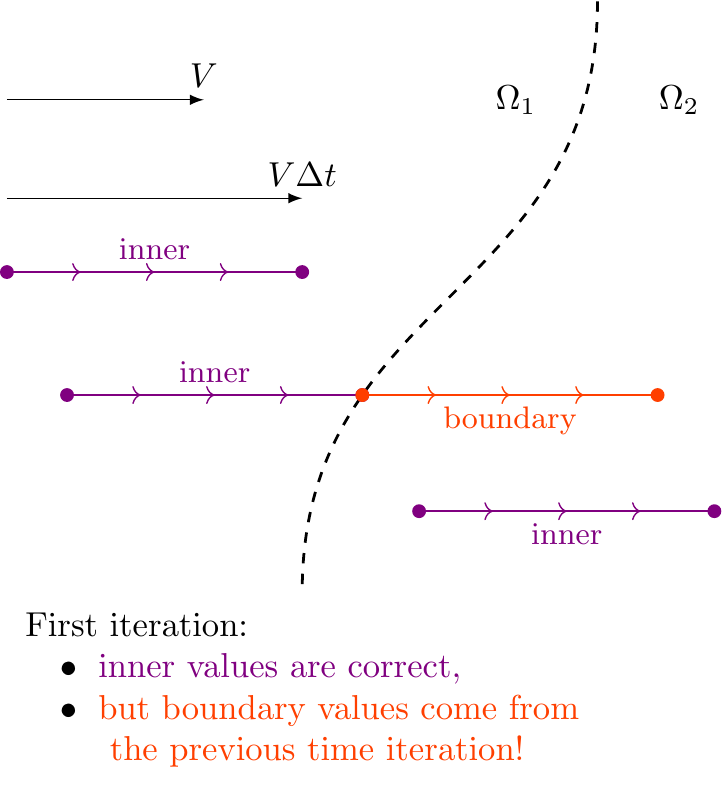}};
		\node (right_graph) [anchor=west, xshift = 1.5cm] at (left_graph.east) {\includegraphics{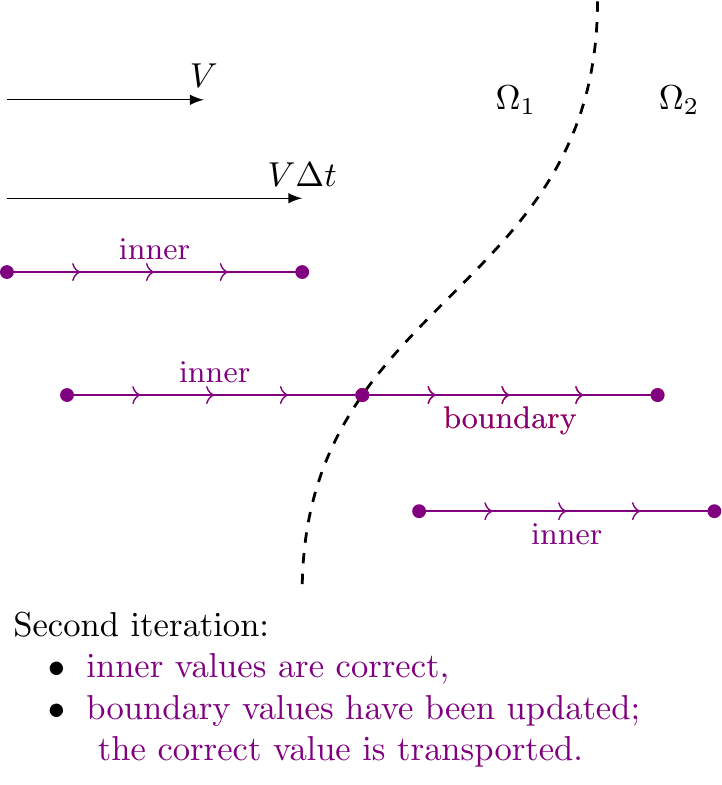}};
		\path (left_graph.south) -- (right_graph.south) node[midway](midway_south){};
		\path (left_graph.north) -- (right_graph.north) node[midway](midway_north){};
		\draw (midway_south) -- (midway_north);
	\end{tikzpicture}

	\caption{%
		Subdomain algorithm, when the subdomain
		decomposition is aligned with the transport velocity.
		In this case, the iterative algorithm reaches the exact solution
		in at most two iterations.
		During the first iteration, in the left panel,
		the boundary values of the subdomains are updated.
		During the second iteration, in the right panel,
		the correct boundary values are transported.
		The purple color corresponds to the transport of a correct value,
		while the red color corresponds to the transport of a wrong value.
	}
	\label{fig:subdomain-two-iter}
\end{figure}

\begin{figure}[!ht]

	\centering

	\makebox[\textwidth][c]{
		\begin{tikzpicture}

			\node (left_graph) at (0, 0) {
				\includegraphics[width=0.38\textwidth]{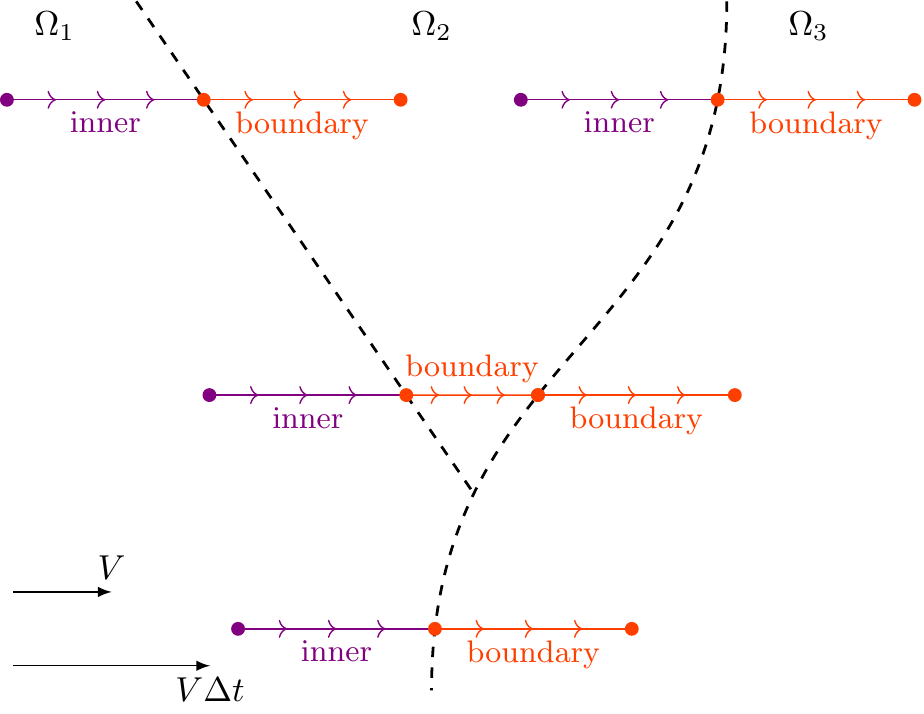}
			};
			\node (center_graph) [anchor=west, xshift = 0cm] at (left_graph.east) {
				\includegraphics[width=0.38\textwidth]{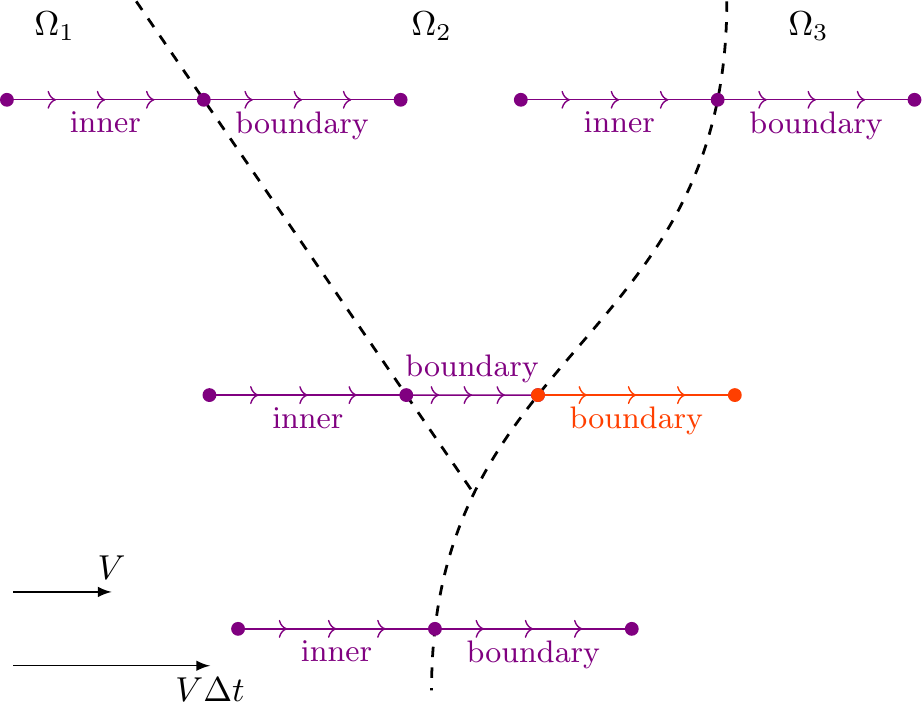}
			};
			\node (right_graph) [anchor=west, xshift = 0cm] at (center_graph.east) {
				\includegraphics[width=0.38\textwidth]{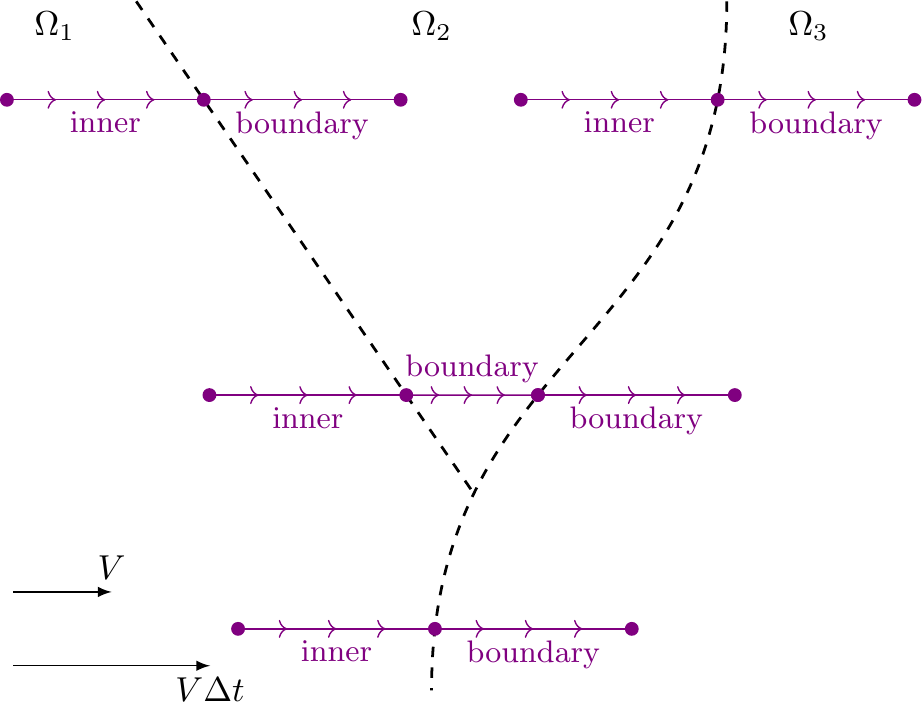}
			};

			\path (left_graph.south) -- (center_graph.south) node[midway](midway_left_south){};
			\path (left_graph.north) -- (center_graph.north) node[midway](midway_left_north){};
			\path (center_graph.south) -- (right_graph.south) node[midway](midway_right_south){};
			\path (center_graph.north) -- (right_graph.north) node[midway](midway_right_north){};
			\draw (midway_left_south) -- (midway_left_north);
			\draw (midway_right_south) -- (midway_right_north);

		\end{tikzpicture}
	}

	\caption{%
		Subdomain algorithm, in a generic subdomain decomposition,
		with corners shared by several subdomains.
		In this case, the iterative algorithm reaches the exact solution
		in at most three iterations.
		First iteration, left panel:
		the boundary value on $\partial\Omega_{2}^{-}$ is updated.
		Second iteration, center panel:
		the boundary value on $\partial\Omega_{3}^{-}$ is updated.
		Third iteration, right panel:
		the correct boundary value is transported.
		The purple color corresponds to the transport of a correct value,
		while the red color corresponds to the transport of a wrong value.
	}

	\label{fig:subdomain-three-iter}
\end{figure}
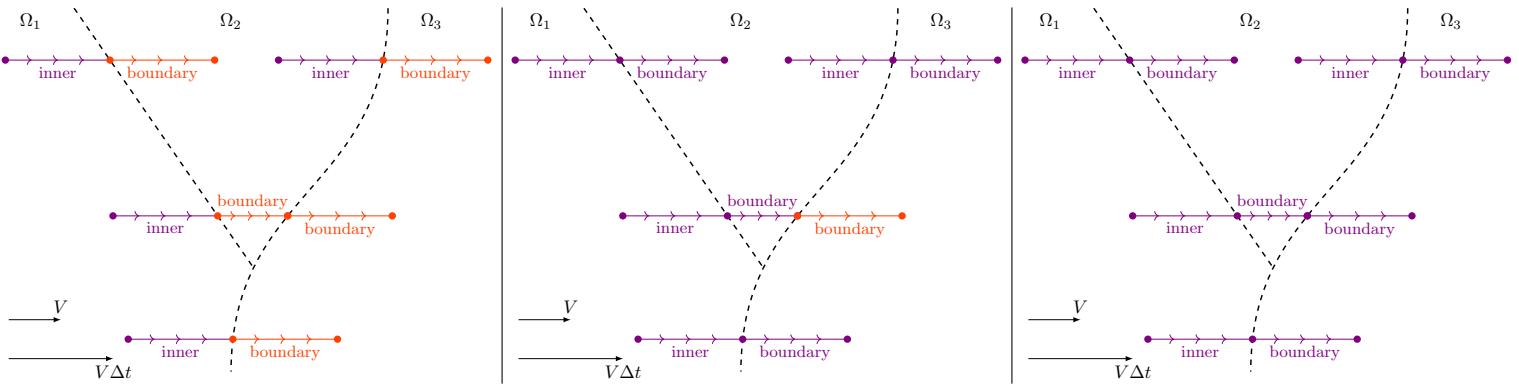

\begin{proof}
	The proof relies on the method of lines.
	It is briefly sketched in \Cref{fig:subdomain-two-iter,fig:subdomain-three-iter}.
\end{proof}

\subsection{Stability}

We have implemented the above iterative algorithm in our Rust code.
The thread-based parallelism within each subdomain is managed, like before, by the
Rust \texttt{rayon} library. Communications between the subdomains
are managed through calls to the MPI (Message Passing Interface) library.
The subdomains are constructed using the METIS graph partitioning tool~\cite{KarKum1998}.

\begin{figure}[!ht]
	\centering
	\includegraphics{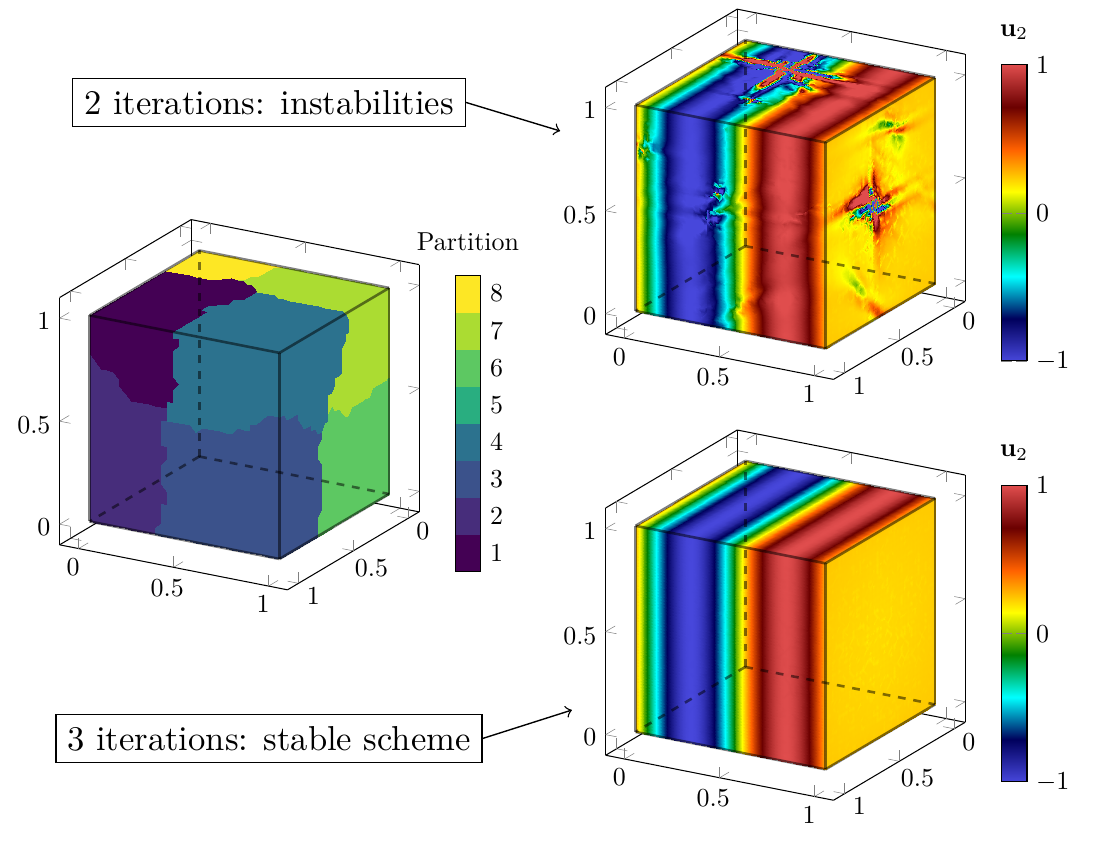}
	\caption{%
		Stability of the subdomain iterative algorithm
		applied on a mesh of the unit cube with 8 subdomains.
		Left: structure of the subdomains;
		Top right: scheme with two iterations;
		Bottom right: scheme with three iterations.
		We observe that the iterative algorithm is stable,
		even with large time steps,
		but that three iterations seem to be necessary.
	}
	\label{fig:iterative-stability}
\end{figure}

First experiments allowed us to verify the stability properties
of the transport solver. They indicate that the number of iterations
of the iterative algorithm is indeed important for the stability
of the method. For a general domain decomposition
and with large time steps, the algorithm is stable provided that three
iterations are performed at each time step. An illustration is
given in \Cref{fig:iterative-stability}.

The objective of the subdomain algorithm was to relax the computational
dependencies and to achieve a better parallel (strong) scaling of
the method. This goal is achieved, as shown in \Cref{tab:mpi-scaling}.
In this table, we compare the time spent in the iterative algorithm
with a varying number of threads and subdomains. We define the efficiency~$e$
of the parallelization as the ratio
between the elapsed time of the algorithm
and the time that we would get with an ideal perfect strong scaling.
The efficiency is perfect if $e=1$.

\begin{table}[!ht]
	\newcolumntype{C}{>{\centering\arraybackslash}p{2.25cm}}
	\centering
	\caption{%
		Multithread and MPI strong scaling test
		on a mesh with about 3.5M elements,
		on a server equipped with an AMD EPYC 7713 x2 (128 physical cores).
		For a computation done with 128 threads,
		we observe that it is better to evenly split the computational work
		between threads and subdomains (green-tinted rows) rather than
		prioritizing threads (red-tinted row).
	}
	\label{tab:mpi-scaling}
	\begin{tabular}{C !{\vrule width -1pt} C !{\vrule width -1pt} C !{\vrule width -1pt} C !{\vrule width -1pt} C}
		\toprule
		\#\,Subdomains & \#\,Threads & \#\,CPU & Time (\SI{}{\second}) & Efficiency \\
		\cmidrule(lr){1-5}
		1              & 1           & 1       & 11350                 & 1.0        \\
		1              & 2           & 2       & 7913.9                & 0.717      \\
		1              & 4           & 4       & 3918.7                & 0.724      \\
		1              & 8           & 8       & 1896.0                & 0.748      \\
		1              & 16          & 16      & 1061.1                & 0.668      \\
		1              & 32          & 32      & 646.12                & 0.549      \\
		1              & 64          & 64      & 424.70                & 0.418      \\
		\rowcolor{badrow!33!white}
		1              & 128         & 128     & 455.70                & 0.195      \\
		2              & 64          & 128     & 250.74                & 0.354      \\
		4              & 32          & 128     & 186.59                & 0.475      \\
		\rowcolor{goodrow!33!white}
		8              & 16          & 128     & 155.03                & 0.572      \\
		\rowcolor{goodrow!33!white}
		16             & 8           & 128     & 155.08                & 0.572      \\
		32             & 4           & 128     & 161.14                & 0.550      \\
		64             & 2           & 128     & 162.38                & 0.546      \\
		128            & 1           & 128     & 162.66                & 0.545      \\
		\bottomrule
	\end{tabular}
\end{table}

We observe, for instance, that with
a single subdomain, the efficiency with $128$ threads drops to $e=0.195$
(as displayed in the red-tinted row),
while with $8$ subdomains and $16$ threads per subdomain, or vice versa,
the efficiency increases to $e=0.572$
(as displayed in the green-tinted rows).
We have thus validated the efficiency of this approach.
Note that the efficiency generally drops with the number of threads,
but this is due to the fact that the mesh remains
too small to provide enough work for each thread.

Of course, the whole algorithm is impacted by a slowdown
imposed by the additional iterations.
However, the weak scaling of the method on a supercomputer
is now certainly ensured for very large computations.
Indeed, explicit subdomain decomposition methods
are known to be well adapted to the architecture of supercomputers,
see for instance \cite{breuer2014sustained,dumbser2018efficient}.

\section{Numerical results}
\label{sec:numerics}

In this section, we present several numerical results obtained with
the kinetic method in three space dimensions.
We apply the method to Maxwell's equations,
and the numerical setup is described in \Cref{sec:setup}.
Several numerical experiments are performed,
namely the propagation of a plane wave in \Cref{sec:plane_wave}
and the simulation of a conductive wire in \Cref{sec:conductive_wire}.
Lastly, a real-world simulation of the interaction of waves emitted by an antenna
with the human body is presented in \Cref{sec:human_body}.

\subsection{Setup of the numerical experiments}
\label{sec:setup}

As a first step, we briefly describe the model used in our numerical experiments,
Maxwell's equations, in \Cref{sec:Maxwell}.
Then, we mention in \Cref{sec:CFL_condition} how the CFL condition is chosen for this 3D problem,
before defining the kinetic velocities $(V_k)_{k \in \{0, \dots, 3\}}$
in \Cref{sec:kinetic_velocities}.
According to \Cref{tab:mpi-scaling}, balancing between
number of subdomains and number of threads leads to the best efficiency.
Unless otherwise mentioned, we use such a setup for each experiment.

\subsubsection{Maxwell's equations}
\label{sec:Maxwell}

Maxwell's equations are a hyperbolic system of conservation laws,
where the vector~$W\in\mathbb{R}^{6}$ of conservative variables is made of
the electric field $E\in\mathbb{R}^{3}$
and the magnetic field $H\in\mathbb{R}^{3}$, as follows:
\[
	W=(E^\intercal,H^\intercal)^\intercal.
\]
Maxwell's equations read
\begin{equation}
	\label{eq:Maxwell_in_PDE_form}
	\begin{dcases}
		\partial_t E - \nabla \times H = - \sigma E, \\
		\partial_t H + \nabla \times E = 0.
	\end{dcases}
\end{equation}
The flux of Maxwell's equations in direction $N \in \mathbb{R}^{3}$ is given by
\[
	Q(W,N)=
	\begin{pmatrix}
		-N\times H \\
		N\times E
	\end{pmatrix},
\]
and we also consider the following source term,
which models a conductive material with conductivity~$\sigma$:
\[
	S(W)=
	\begin{pmatrix}
		- \sigma E \\
		0
	\end{pmatrix}.
\]

\subsubsection{CFL condition}
\label{sec:CFL_condition}

In order to properly compare methods, we have to define the CFL condition.
The reader is referred to \cite{gerhard2022unconditionally}
for a more in-depth discussion on CFL conditions for DG methods.
Here, we define the time step $\Delta t$ as follows:
\begin{equation}
	\label{eq:def_CFL}
	\Delta t = \beta \frac {1}{\lambda_{\text{max}}} h_{\text{min}}.
\end{equation}
In this definition, $\lambda_{\text{max}}$ is the maximum eigenvalue of the Jacobian matrix
of the flux associated to Maxwell's equations (here, $\lambda_{\text{max}} = 1$).
In addition, $h_{\text{min}}$ is defined as the size of the smallest cell in the mesh:
\begin{equation*}
	h_{\text{min}} = \min_{L \in \mathcal{M}} \text{size}(L)
	\text{, \; where \quad}
	\text{size}(L) = \frac{\text{volume}(L)}{\text{surface}(\partial L)}.
\end{equation*}
Finally, $\beta$ is the CFL number.
The maximum possible value for $\beta$ is constrained by the scheme under consideration;
for classical explicit DG schemes, $\beta$ must be of the order of $1$ to get stability.
As we will see in the numerical experiments, we are able to take $\beta$ as large as we want
without loss of stability (but incurring a loss in precision).
This hold whatever value of $\omega < 2$ is chosen,
so we take $\omega = 2 - 10^{-12}$ in order to ensure stability while remaining second-order accurate.

\subsubsection{Kinetic velocities}
\label{sec:kinetic_velocities}

The last ingredient needed to define the scheme is the set of kinetic velocities.
The simplest choice is to choose the following velocities, called ``D3Q4'' in the Lattice-Boltzmann community:
\begin{equation*}
	V_0
	=
	\begin{pmatrix}
		\lambda \\
		\lambda \\
		\lambda
	\end{pmatrix},
	\quad
	V_1
	=
	\begin{pmatrix}
		\lambda   \\
		- \lambda \\
		- \lambda
	\end{pmatrix},
	\quad
	V_2
	=
	\begin{pmatrix}
		- \lambda \\
		\lambda   \\
		- \lambda
	\end{pmatrix},
	\quad
	V_3
	=
	\begin{pmatrix}
		- \lambda \\
		- \lambda \\
		\lambda
	\end{pmatrix},
\end{equation*}
with $\lambda = \sqrt{3}$ to satisfy the subcharacteristic condition.
Note that, with this velocity set, the equilibrium functions $M_k(W)$ read,
for all $k \in \{0, 1, 2, 3\}$:
\begin{equation*}
	M_k(W) = \frac W 4 + \frac{Q(W, V_k)}{4 \lambda^2}.
\end{equation*}


\subsection{Plane wave}
\label{sec:plane_wave}

\begin{figure}[!ht]%
	\makebox[\textwidth][c]{
		\includegraphics[width=0.3\textwidth]{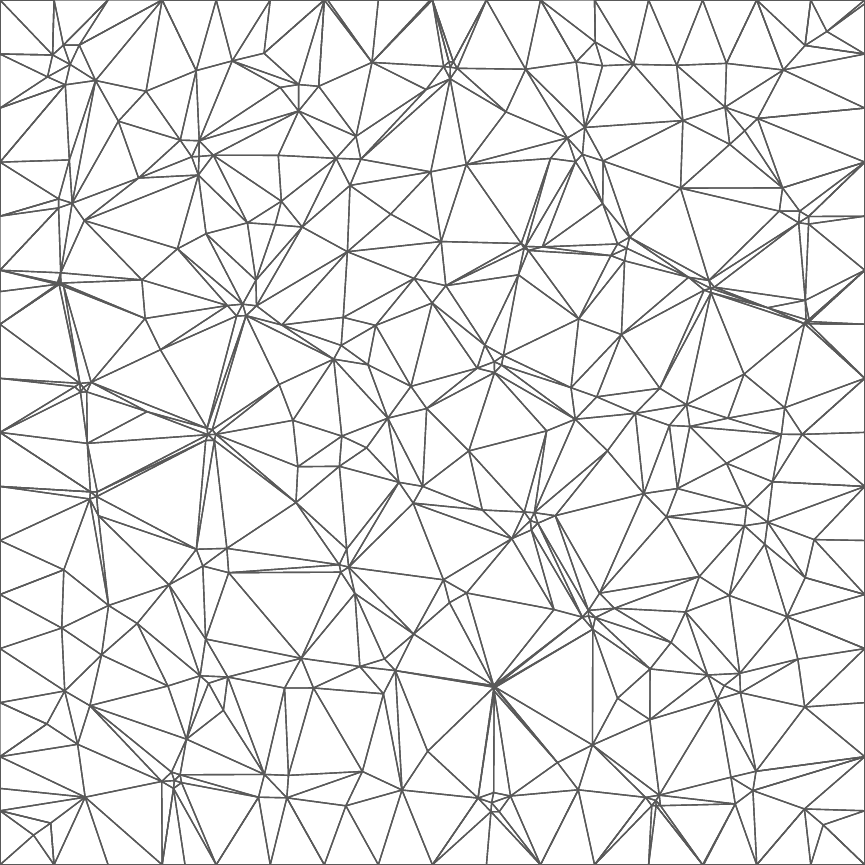}
		\quad
		\includegraphics[width=0.3\textwidth]{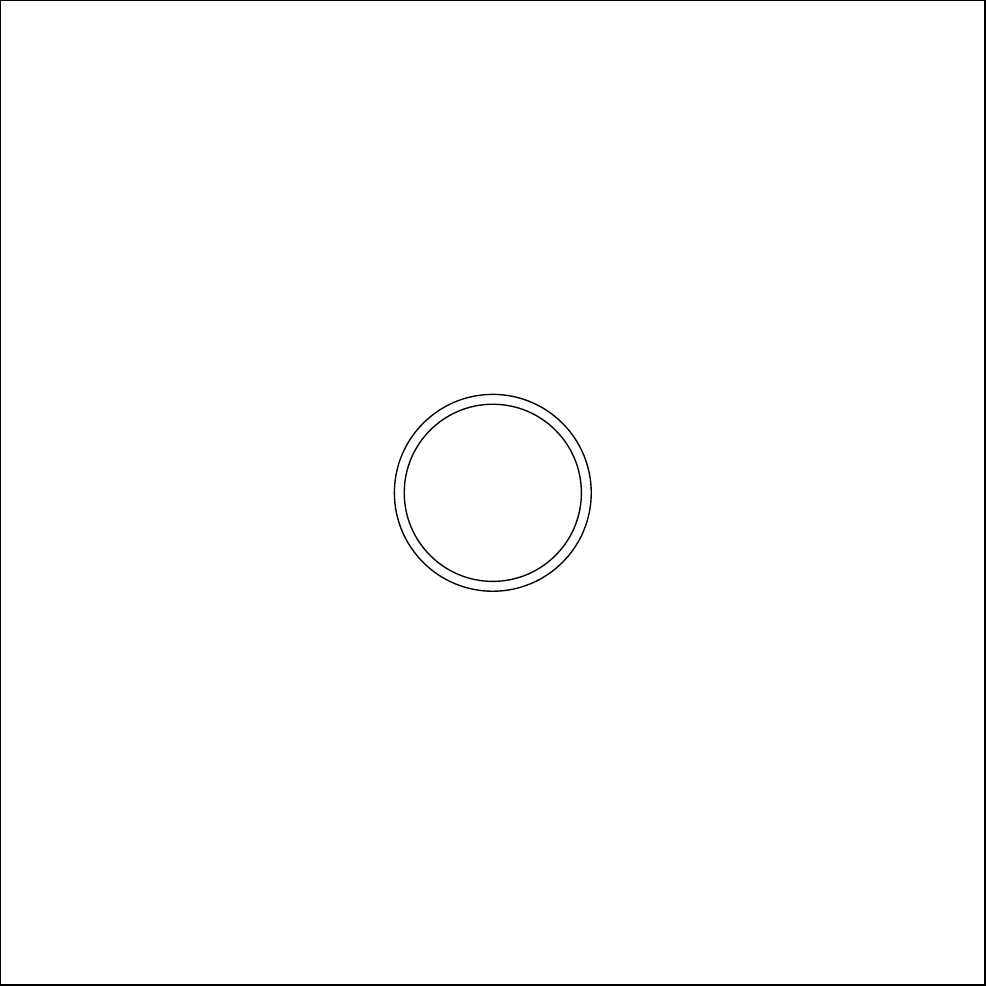}
		\quad
		\includegraphics[width=0.3\textwidth]{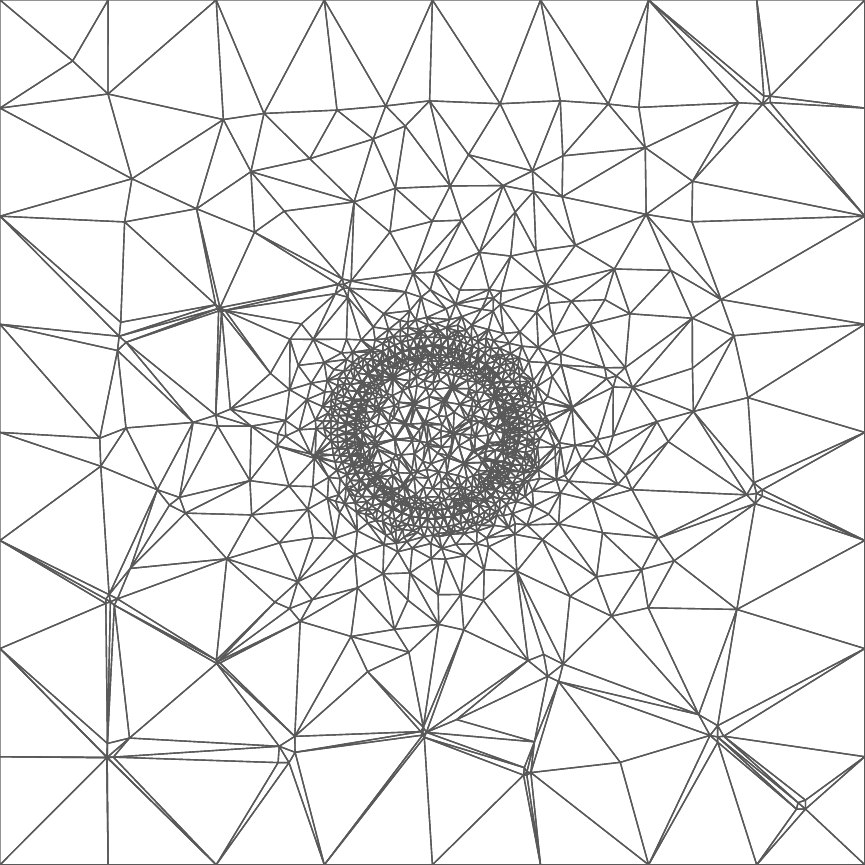}
	}
	\caption{%
		Two meshes of the unit cube, sliced at $x_2 = 0.5$.
		The left panel contains a mesh with uniformly spaced tetrahedra,
		labeled $\mathcal{M}_1$,
		while the other two focus on a mesh with local refinement at the center,
		labeled $\mathcal{M}_2$.
		The center panel contains the geometry of this local refinement
		(in the shape of a torus),
		while the right panel depicts the locally refined mesh itself.
	}
	\label{fig:Mesh}
\end{figure}

In order to validate the subdomain decomposition,
we first run an experiment already performed in \cite{gerhard2022unconditionally}
without the subdomain decomposition.
We expect the results to be almost the same.
The computational domain is the unit cube.
For this test, the conductivity $\sigma$ is set to zero.
We consider two meshes, represented on \Cref{fig:Mesh},
one with a uniform cell size (labeled $\mathcal{M}_1$)
and one with a non-uniform cell size (labeled $\mathcal{M}_2$).
We compute, with the meshes described above,
the propagation of a plane wave with frequency~$\nu$.
The exact solution therefore is
\begin{equation*}
	W(X, t) =
	\begin{pmatrix}
		0                           \\
		0                           \\
		\cos(2 \pi \nu (x_1 - t))   \\
		0                           \\
		- \cos(2 \pi \nu (x_1 - t)) \\
		0
	\end{pmatrix},
\end{equation*}
and we prescribe this exact solution on the boundary with Dirichlet boundary conditions.

We check the CFL-less feature in \Cref{tab:cfl_less_cube_16}.
To that end, we define the error $e_r(\nu)$ between the exact and approximate solutions.
For a more precise description of this error and of the whole setup,
the reader is referred to~\cite{gerhard2022unconditionally}.
In this table, we verify that the scheme is stable,
even at very high CFL numbers $\beta$, for both meshes.
In addition, for fixed $\Delta t$
(which is not the same as fixed $\beta$ according to the definition \eqref{eq:def_CFL} of $\Delta t$),
the scheme is about as precise on both meshes.
Note that, for a standard third order explicit DG scheme to be stable in this configuration,
a CFL condition $\beta \leq 1.85$ is required.
However, the kinetic scheme remains stable and precise when using larger time steps,
for which an explicit DG scheme would not be stable any longer.
For instance, on the locally refined mesh $\mathcal{M}_2$,
we are able to take a CFL number $\beta = 37 = 20 \times 1.85$
while retaining about the same error as with $\beta = 1.85$.
Moreover, the results are the same as
the single-subdomain version from~\cite{gerhard2022unconditionally},
which further validates our approach.

\begin{table}[!ht]
	\caption{%
		Numerical results for the experiment described
		in \cref{sec:plane_wave} and pictured in \Cref{fig:Mesh}.
		For both meshes, we collect the time step $\Delta t$ and the error $e_r$
		with respect to the value of the CFL number and to the choice of frequency.
	}
	\label{tab:cfl_less_cube_16}
	\newcolumntype{C}{>{\centering\arraybackslash}p{1.5cm}}
\newcolumntype{R}{>{\raggedleft\arraybackslash}p{1.5cm}}
\newcolumntype{L}{>{\raggedright\arraybackslash}p{1.5cm}}

\pgfplotstableread{./data/article_cfl_less_cube_16_maxwell_cos_t1.dat}\datauniformmesh
\pgfplotstableread{./data/article_cfl_less_cube_8_torus_128_maxwell_cos_t1.dat}\datanonuniformmesh

\pgfplotstablecreatecol[copy column from table={\datanonuniformmesh}{[index] 2}] {dt_non_uniform} {\datauniformmesh}
\pgfplotstablecreatecol[copy column from table={\datanonuniformmesh}{[index] 3}] {err_cos2_non_uniform} {\datauniformmesh}
\pgfplotstablecreatecol[copy column from table={\datanonuniformmesh}{[index] 4}] {err_cos5_non_uniform} {\datauniformmesh}

\makebox[\textwidth][c]{
	\small
	\pgfplotstabletypeset[
		empty cells with={---},
		columns={cfl, dt, err_cos2, err_cos5, dt_non_uniform, err_cos2_non_uniform, err_cos5_non_uniform},
		columns/cfl/.style   = {fixed, fixed zerofill, precision = 2, column name = CFL,        column type = {C}},
		columns/dt/.style    = {fixed, fixed zerofill, precision = 5, column name = $\Delta t$, column type = {C}},
		columns/dt_non_uniform/.style    = {fixed, fixed zerofill, precision = 5, column name = $\Delta t$, column type = {C}},
		columns/err_cos2/.style   = {fixed, fixed zerofill, precision = 5, column name = $\nu = 2$,      column type = {C}},
		columns/err_cos2_non_uniform/.style   = {fixed, fixed zerofill, precision = 5, column name = $\nu = 2$,      column type = {C}},
		columns/err_cos5/.style   = {fixed, fixed zerofill, precision = 5, column name = $\nu = 5$,      column type = {C}},
		columns/err_cos5_non_uniform/.style   = {fixed, fixed zerofill, precision = 5, column name = $\nu = 5$,      column type = {C}},
		columns/CPU_1/.style   = {fixed, fixed zerofill, precision = 2, column name = CPU ($s$),  column type = {C}},
		columns/CPU_24/.style   = {fixed, fixed zerofill, precision = 2, column name = CPU ($s$),  column type = {C}},
		every head row/.style={output empty row, before row={\toprule%
						&  \multicolumn{3}{c}{mesh $\mathcal{M}_1$}  & \multicolumn{3}{c}{mesh $\mathcal{M}_2$} \\
						\cmidrule(lr){2-4}\cmidrule(lr){5-7}}},
		every head row/.append style={after row={CFL $\beta$ &  $\Delta t$ & $e_r(2)$ & $e_r(5)$ & $\Delta t$ & $e_r(2)$ & $e_r(5)$ \\
						\cmidrule(lr){1-1}%
						\cmidrule(lr){2-2}%
						\cmidrule(lr){3-4}%
						\cmidrule(lr){5-5}%
						\cmidrule(lr){6-7}}},
	]{\datauniformmesh}
}

\end{table}

\begin{figure}[!ht]%
	\centering
	\includegraphics[width=\textwidth]{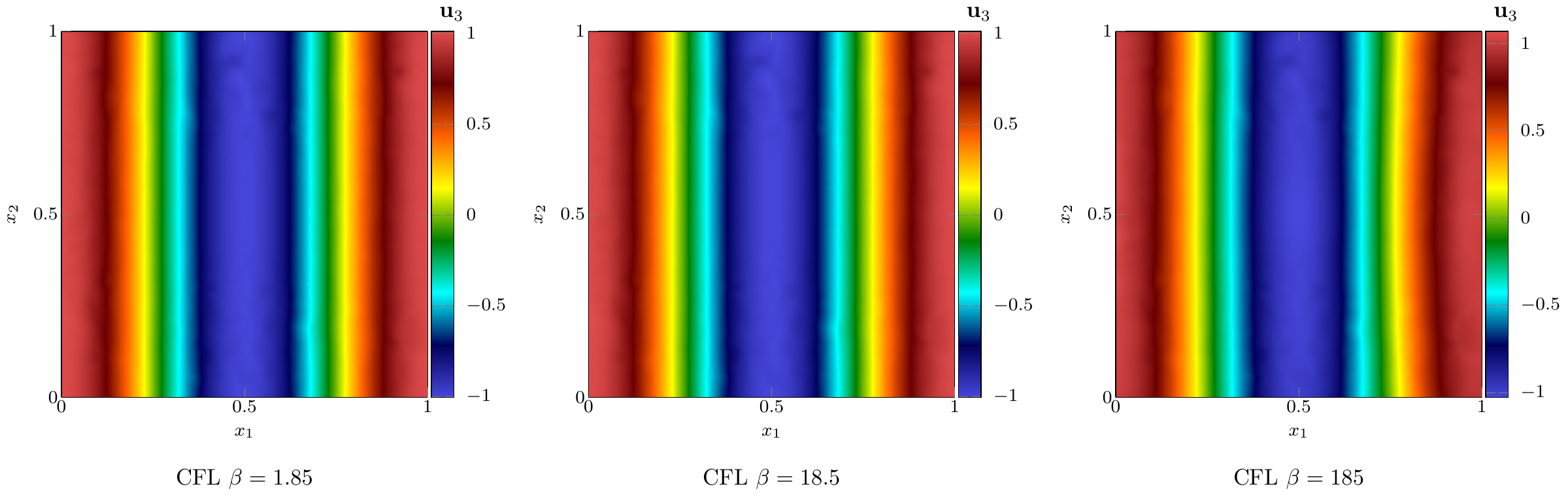}
	\includegraphics[width=\textwidth]{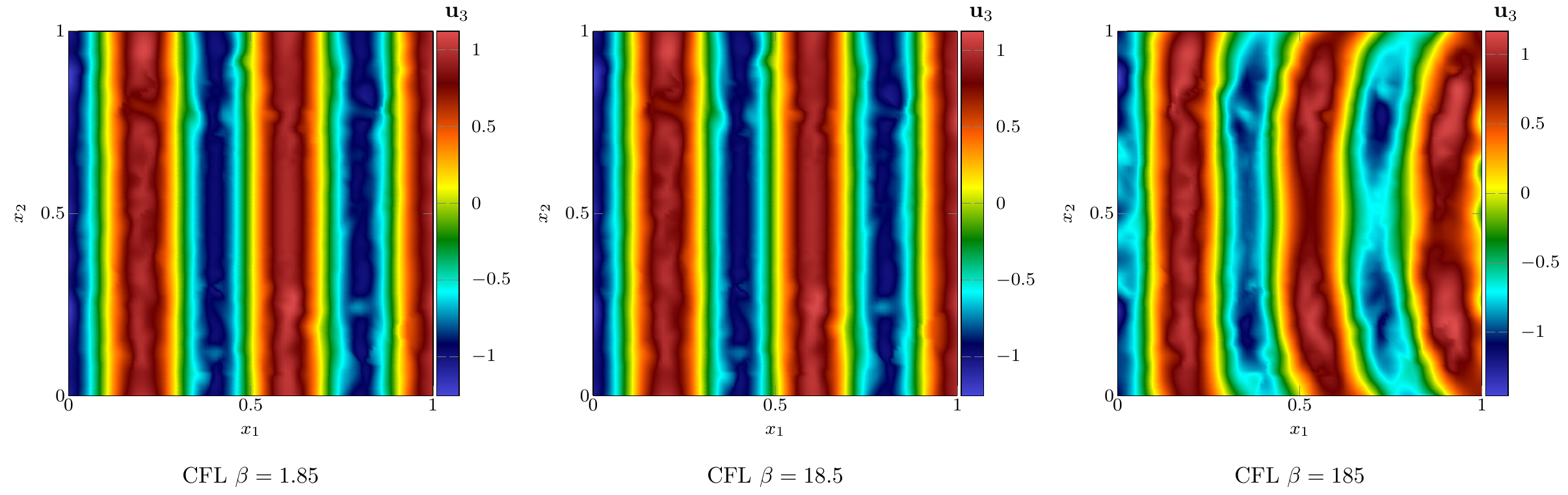}
	\caption{%
		Plane wave propagation from \Cref{sec:plane_wave}
		on the locally refined mesh $\mathcal{M}_2$:
		depiction, for $x_3=0.5$, of the third component $W_3(X,t) = E_3(X,t)$
		of the approximate solution.
		We simulated plane waves with frequencies
		$\nu = 2$ (top panels) and $\nu = 5$ (bottom panels).
		From left to right, we have set the CFL number $\beta$ to $1.85$, $18.5$ and $185$.
	}
	\label{fig:cfl_less_cube_8_torus_128}
\end{figure}

Color plots of the numerical results are given in
\Cref{fig:cfl_less_cube_8_torus_128}.
These plots illustrate the stability of the computations at high CFL numbers,
and we note that the approximate solution remains stable
even for extremely large values of the CFL number $\beta$.
Of course, the accuracy of the computation depends
on the frequency of the plane wave: the more the solution oscillates,
the smaller the time step should be in order to accurately capture the oscillations.

\subsection{Conductive wire}
\label{sec:conductive_wire}

\begin{figure}[!ht]
	\centering
	\includegraphics[width=0.45\textwidth]{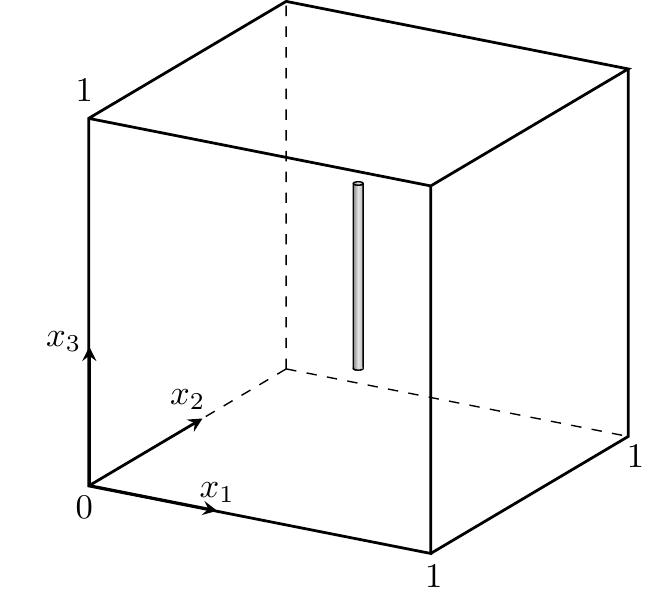}
	\qquad
	\includegraphics[width=0.405\textwidth]{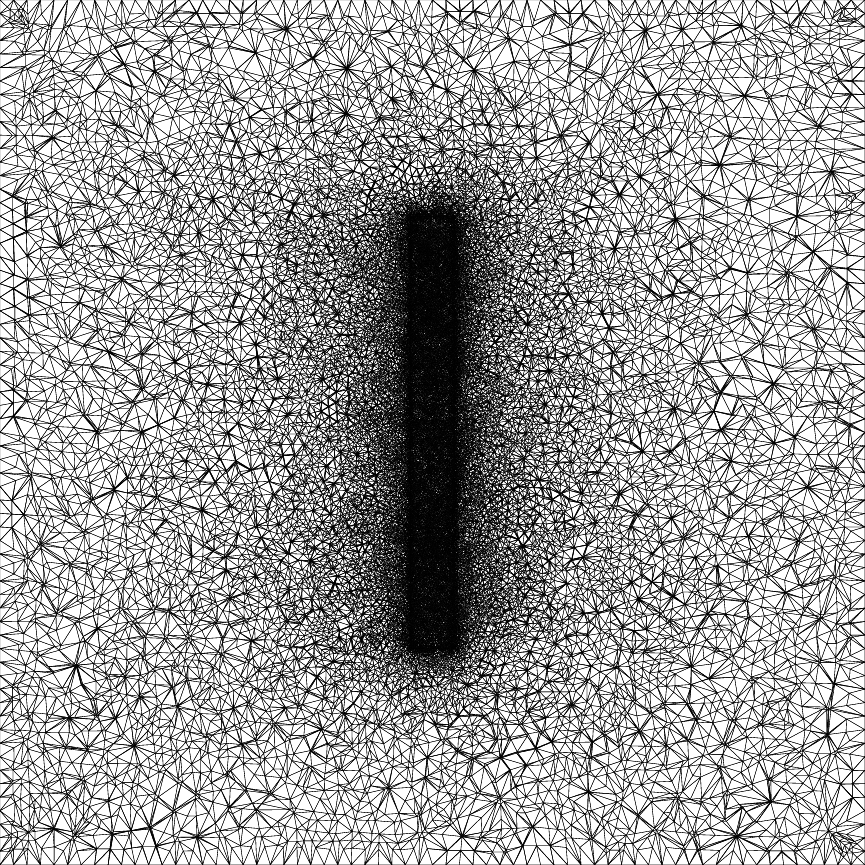}
	\caption{%
		Geometry (left panel) and mesh (right panel)
		for the conductive wire test from
		\Cref{sec:conductive_wire}.
		The mesh is locally refined around the wire,
		where the conductivity $\sigma$ is nonzero.
	}
	\label{fig:wire-geom}
\end{figure}

In this test, we activate the source term,
i.e., we take a nonzero $\sigma$ (at least in some part of the domain).
We consider a small electric wire located in
the middle of the computational domain.
The unstructured mesh of the unit cube conforms with the small wire,
which means large cells far from the wire
and locally refined cells close to the wire (see \Cref{fig:wire-geom}).
Once again, this experiment was also performed in \cite{gerhard2022unconditionally},
and we present it here to further validate the transport algorithm on the subdomain decomposition.

In order to validate the proposed methodology,
we compare our DG method with a well-validated FDTD
(Finite-Difference Time-Domain) solver~\cite{guiffaut2011new}, based on
the Yee scheme~\cite{yee1966numerical}, that can handle electric wires.
This FDTD solver requires a uniform Cartesian grid,
which means it is easy to parallelize.
Therefore, the fine mesh within the antenna implies
a uniformly fine mesh everywhere in the domain;
in practice, we use $1000^3 = 10^9$ cells.

We also compare our results to an explicit RK2-DG solver,
called CLAC (Computation Laws on mAny Cores).
This solver is contained within a well-validated code,
parallelized on a GPU, as opposed to our implicit kinetic solver.

The differences between the solvers are summarized in
\Cref{tab:solvers_summary}.

\begin{table}[!ht]
	\caption{%
		Summary of the algorithms used in this section:
		algorithm type, programming language and parallelization type.
	}
	\label{tab:solvers_summary}
	\centering
	\begin{tabular}{cccc}
		\toprule
		Solver   & Algorithm           & Language                      & Parallelization            \\ \cmidrule(lr){1-4}
		FDTD     & finite differences  & \texttt{Fortran}              & CPU, distributed           \\
		CLAC     & explicit RK2-DG     & \texttt{C++}, \texttt{OpenCL} & GPU                        \\
		KOUGLOFV & implicit kinetic DG & \texttt{Rust}                 & CPU, shared \& distributed \\ \bottomrule
	\end{tabular}
\end{table}

To set up the numerical experiment, a plane wave pulse is sent through the vacuum, towards the wire.
This amounts to solving Maxwell's equations
with the conductivity source term $S(W)=(\sigma E,0)^\intercal$,
where the conductivity $\sigma$ vanishes outside the wire.
To define the initial and boundary conditions,
we consider the following exact solution of Maxwell's equations without source term:
\begin{equation*}
	W(X, t) =
	\begin{pmatrix}
		0                     \\
		0                     \\
		\psi(x_2 - x_c - t)   \\
		- \psi(x_2 - x_c - t) \\
		0                     \\
		0
	\end{pmatrix},
\end{equation*}
where $x_c = 0.25$ and where $\psi$ is a compactly supported bump function:
\begin{equation*}
	\psi(X) =
	\begin{dcases}
		\exp \left( 1 - \frac 1 {1 - \frac { \|X\| } \eta} \right) & \text{ if } \|X\| < \eta, \\
		0                                                          & \text{ otherwise,}
	\end{dcases}
\end{equation*}
with $\eta = 0.25$ the size of the bump.
Then, the initial condition is $W(X, 0)$,
the boundary conditions consist in imposing the solution $W(X, t)$ at the boundaries.
Note that $W(X, t)$ is not an exact solution of the problem with source term;
however, the antenna is far enough away from the boundaries for this fact not to matter when running simulations.

Note that, in this case, the application of the source term in the scheme,
according to \eqref{eq:source_term_application_to_W}, reads:
\begin{equation*}
	\label{eq:source_term_application_to_W_conductivity}
	\left\{
	\begin{aligned}
		\frac{E(\cdot,\Delta t^{+})-E(\cdot,\Delta t^{-})}{\Delta t} & =\frac{\sigma E(\cdot,\Delta t^{+})+\sigma E(\cdot,\Delta t^{-})}{2}, \\
		\frac{H(\cdot,\Delta t^{+})-H(\cdot,\Delta t^{-})}{\Delta t} & =0.
	\end{aligned}
	\right.
\end{equation*}
This is a linear equation in the unknown $E(\cdot,\Delta t^{+})$,
which gives:
\begin{equation}
	\label{eq:W_after_conductivity_source_term}
	\begin{dcases}
		E(\cdot,\Delta t^{+}) = \mu E(\cdot,\Delta t^{-}), \\
		H(\cdot,\Delta t^{+}) = H(\cdot,\Delta t^{-}),
	\end{dcases}
	\text{\quad where \quad}
	\mu = \dfrac{1 - \sigma \dfrac{\Delta t}{2}}{1 + \sigma \dfrac{\Delta t}{2}}.
\end{equation}

We now apply the solvers to two different cases.
Within the wire, we consider two values of the conductivity~$\sigma$:
first, a small conductivity $\sigma=3$ in \Cref{sec:conductive_wire_sigma_3};
then, an infinite conductivity $\sigma \to +\infty$
in \Cref{sec:conductive_wire_sigma_inf}.
For the KOUGLOFV solver, we take a CFL number $\beta = 7$.
Recall that the scheme is stable whatever the value of~$\beta$,
but this choice yields a good compromise
between precision and computation speed.
Since the three solvers give comparable results,
a computation time comparison is proposed in \Cref{sec:conductive_wire_CPU_time}.

\subsubsection{Wire with a low conductivity}
\label{sec:conductive_wire_sigma_3}

\begin{figure}[!ht]
	\centering
	\includegraphics[width=\textwidth]{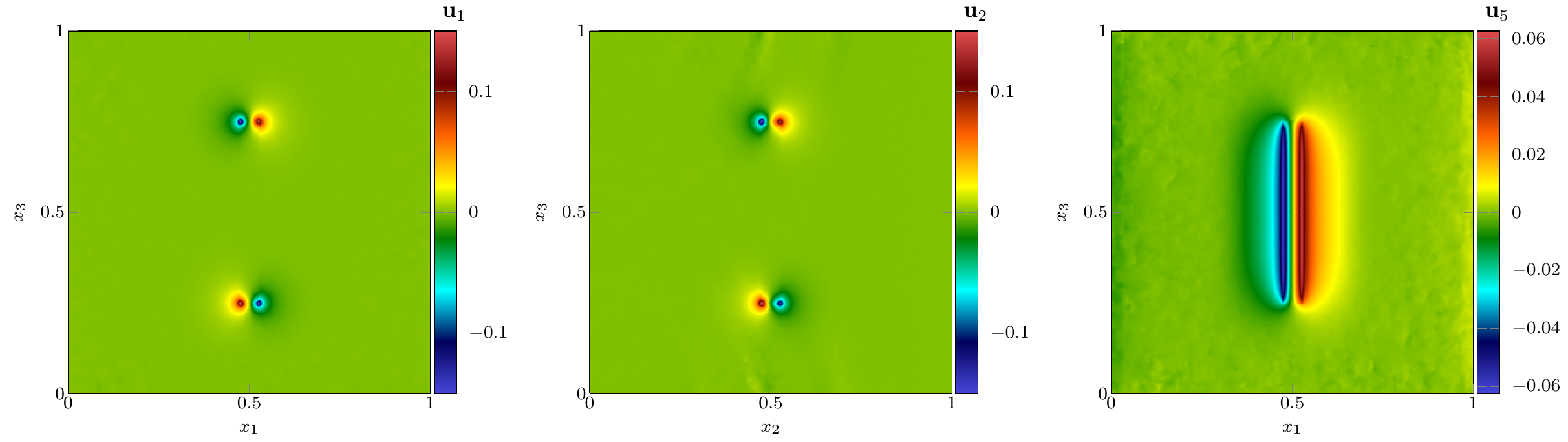}
	\caption{%
		Conductive wire simulation from \Cref{sec:conductive_wire_sigma_3}, with $\sigma = 3$, solution at $t = 0.75$:
		left panel: $\restr{E_1}{x_2 = 0.5}$;
		center panel: $\restr{E_2}{x_1 = 0.5}$;
		right panel: $\restr{H_2}{x_2 = 0.5}$.
	}
	\label{fig:wire-results}
\end{figure}

We first consider $\sigma = 3$.
We display on \Cref{fig:wire-results}
the approximate solution obtained by the KOUGLOFV solver.
We observe a good agreement with the expected results,
since the electric charge is concentrated at the ends of the antenna,
and the magnetic field rotates around the antenna.

\begin{figure}[!ht]
	\centering
	\includegraphics{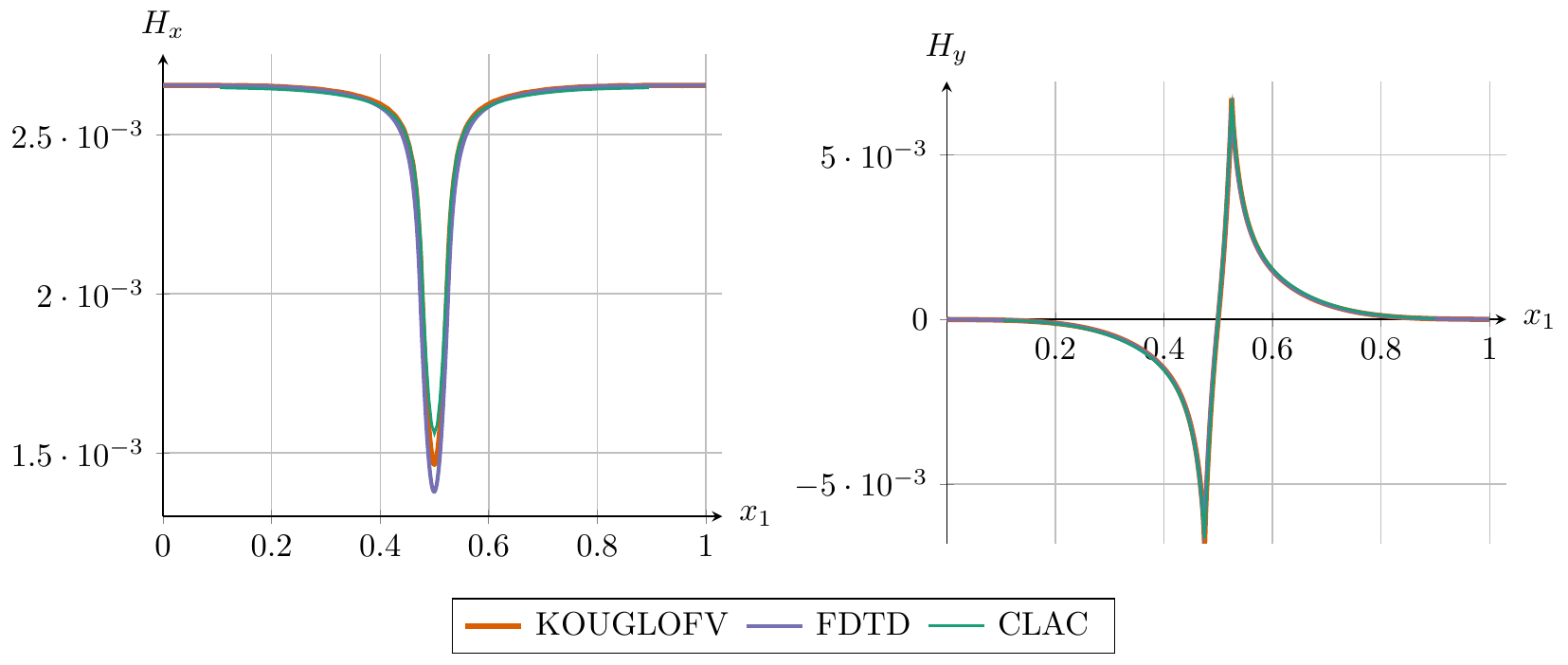}
	\caption{%
		Conducting wire with small conductivity $\sigma = 3$,
		comparison between the implicit DG solver (KOUGLOFV),
		the FDTD solver and a standard explicit DG solver (CLAC).
	}
	\label{fig:wire-sigma3}
\end{figure}

Then, to compare the KOUGLOFV solver with the FDTD and CLAC solvers,
the first and second components of the magnetic field
$H_{x}$ and $H_{y}$, along the line $y=z=1/2$ and when the pulse reaches
the wire, are represented on \Cref{fig:wire-sigma3}.
Note that, for CLAC, the results are displayed on the domain $(0.1, 0.9)$
instead of $(0, 1)$.
This is due to a technical limitation of the CLAC code related to the boundary conditions,
which is not present for the FDTD and KOUGLOFV solvers.

We observe an excellent agreement with the FDTD and CLAC solvers.
Let us mention that, to obtain these results,
our implicit kinetic DG solver does not apply any charge conservation
correction, while this is generally considered to be necessary for
such simulations,
see for instance~\cite{munz2000divergence,crestetto2012resolution}.
Further investigation is needed for understanding this good behavior,
which is perhaps linked to the test case under consideration.

\subsubsection{Infinitely conductive wire}
\label{sec:conductive_wire_sigma_inf}

\begin{figure}[!ht]
	\includegraphics[width=\textwidth]{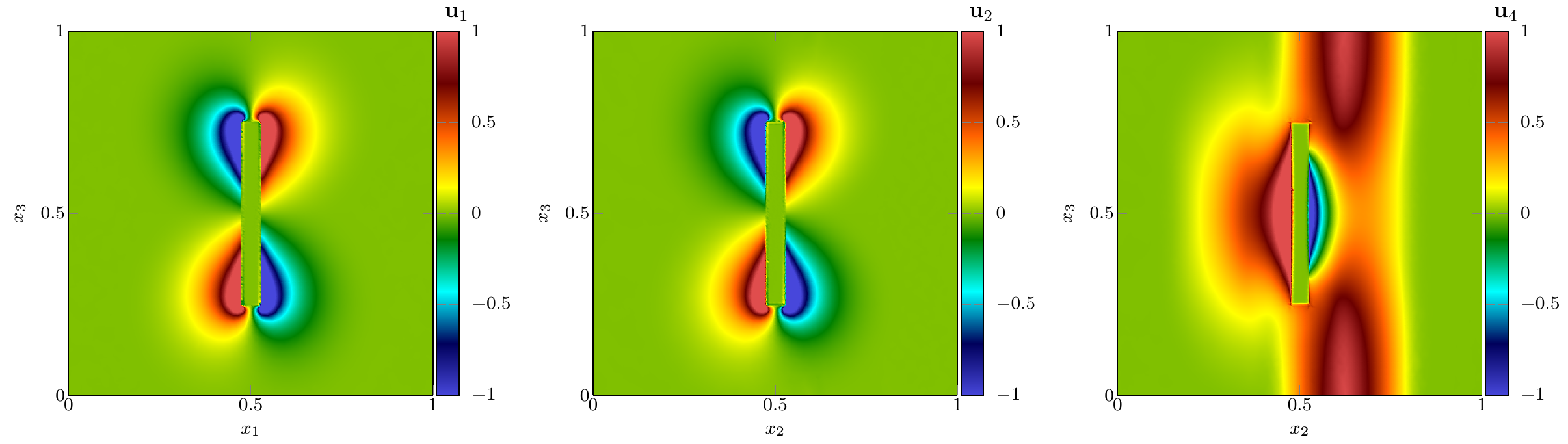}
	\caption{%
		Infinitely conductive wire simulation from \Cref{sec:conductive_wire_sigma_inf},
		with $\sigma \to +\infty$,
		solution at $t = 0.375$:
		left panel: $\restr{E_1}{x_2 = 0.5}$;
		center panel: $\restr{E_2}{x_1 = 0.5}$;
		right panel: $\restr{H_1}{x_1 = 0.5}$.
	}
	\label{fig:numerical_source_simulation_stiff}
\end{figure}

In this section,
we consider a test with a large conductivity $\sigma \to +\infty$.
In this case, the source term is very stiff,
but is handled without issue by our method,
thanks to the implicit source term treatment.
Indeed, recall equation \eqref{eq:W_after_conductivity_source_term}:
when $\sigma \to +\infty$, we get
\begin{equation*}
	\label{eq:W_after_infinite_conductivity_source_term}
	\begin{dcases}
		E(\cdot,\Delta t^{+}) = - E(\cdot,\Delta t^{-}), \\
		H(\cdot,\Delta t^{+}) = H(\cdot,\Delta t^{-}).
	\end{dcases}
\end{equation*}
This behavior is consistent with what would happen in
perfect electrical conductors (PECs).
In practice, we take $\sigma = 10^{12}$ to mimic infinity.
This gives the same results as writing the source term in the $\sigma = +\infty$ limit.
The numerical results are depicted in \Cref{fig:numerical_source_simulation_stiff}.
We observe a good agreement with the expected values of the electric and magnetic fields.
Indeed, the electric field (left and center panels) remains constant, equal to zero, within the wire.
In addition, the right panel clearly shows that the wave
(traveling from left to right) has just passed the wire,
but has not created a magnetic field within the wire.

\begin{figure}[!ht]
	\centering
	\includegraphics{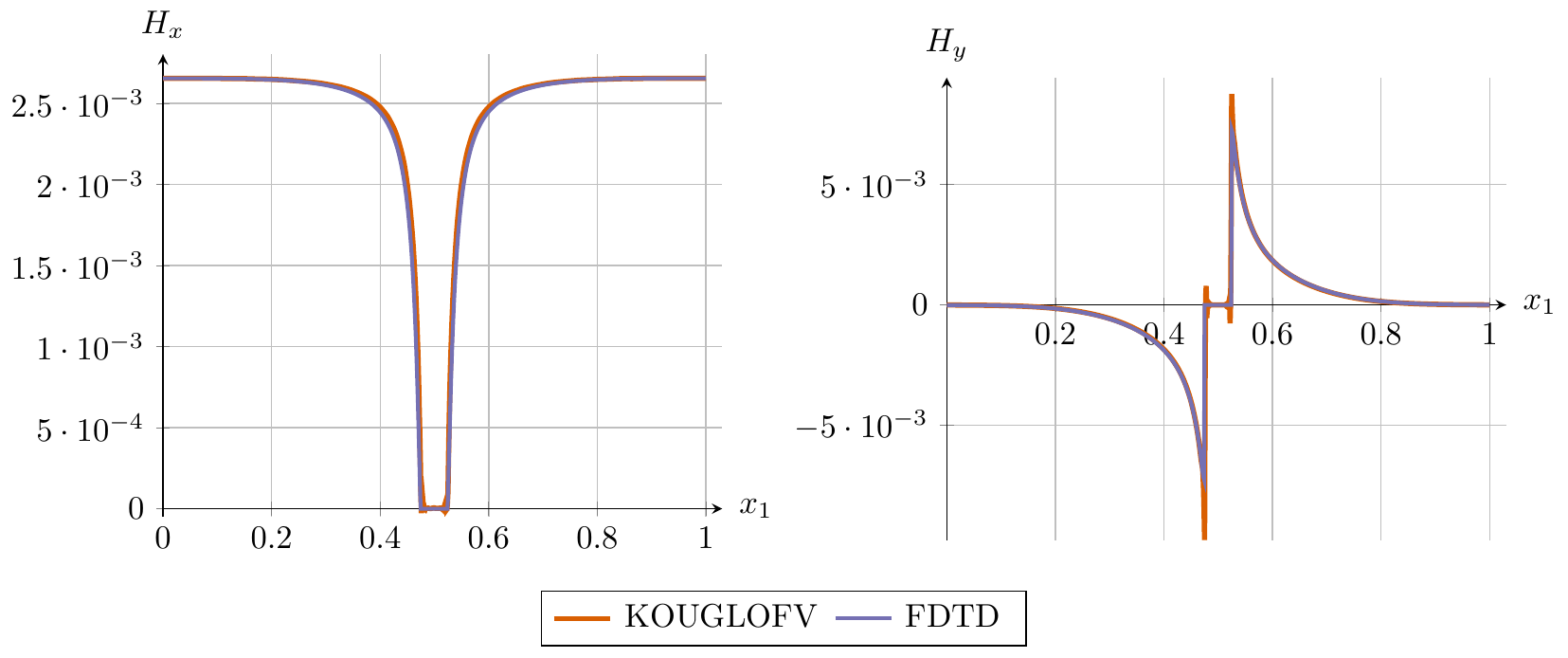}
	\caption{%
		Conducting wire with infinite conductivity $\sigma \to +\infty$,
		comparison between the implicit DG solver (KOUGLOFV),
		and the FDTD solver.
	}
	\label{fig:wire-sigma-inf}
\end{figure}

Contrary to this, in the explicit methods,
a vanishingly small CFL condition (proportional to $\frac 1 \sigma$)
is needed for stability.
This is, of course, unusable in practice,
and a specific PEC formulation has to be used.
This is another advantage of our method:
there is no need to implement complex PEC conditions
to handle infinite conductivities.
We compare the approaches on \Cref{fig:wire-sigma-inf},
where we observe very good agreement between the solutions.
The implicit DG solver presents a few small oscillations close to the wire.
They can mostly be attributed to visualization,
since the software used to create the slices (Paraview, see \cite{Aya2015})
performs an interpolation.
A large part of the oscillation amplitude is due to these visualization artifacts,
and a small part is due to standard dispersion effects
(since we consider a high-order scheme with a large CFL number).


\subsubsection{Computation time comparison}
\label{sec:conductive_wire_CPU_time}

Finally, in light of the similar results obtained in both the low and the infinite conductivity cases,
we report, in \Cref{tab:FDTD_vs_KOUGLOFV_CPU_time}, the computation time taken by each solver.
For the FDTD solver, the sheer size of the mesh makes
obtaining the approximate solution quite a bit slower
than with the KOUGLOFV solver.
The CLAC solver is much faster than the FDTD solver,
but because of the restrictive CFL condition,
obtaining the solution is still faster with the KOUGLOFV solver.
Note that the CLAC solver is parallelized on a GPU,
while the KOUGLOFV solver only involves a multi-core CPU.
Even so, the KOUGLOFV solver manages to be significantly faster than the CLAC solver,
and it takes a comparable computation time when ran on fairly ancient desktop CPUs.

\begin{table}[!ht]
	\caption{%
		CPU time taken to solve the conductive wire test case:
		comparison between the FDTD method,
		the explicit RK2-DG solver (CLAC),
		and the kinetic DG solver (KOUGLOFV).
		For KOUGLOFV, we take a CFL number $\beta = 7$.
	}
	\label{tab:FDTD_vs_KOUGLOFV_CPU_time}
	\centering
	\begin{tabular}{cccc}
		\toprule
		solver                    & number of cells       & hardware platform                                      & computation time \\ \cmidrule(lr){1-4}
		FDTD                      & $1000^3$              & AMD EPYC 7302$\times$2, 32 cores, \SI{3}{\giga\hertz}  & 15.3 hours       \\[\medskipamount]
		CLAC                      & 236k                  & Nvidia GeForce GTX 1070
		                          & 5.98 minutes                                                                                      \\[\medskipamount]
		\multirow{2}{*}{KOUGLOFV} & \multirow{2}{*}{516k} & Intel i7-5820K, 6 cores, \SI{3.3}{\giga\hertz}         & 17.1 minutes     \\
		                          &                       & AMD EPYC 7713$\times$2, 128 cores, \SI{2}{\giga\hertz} & 78.7 seconds     \\ \bottomrule
	\end{tabular}
\end{table}

\subsection{Real-world simulation: interaction of waves from an antenna with an anthropomorphic mannequin}
\label{sec:human_body}

In this configuration, we compare our DG method
to the CLAC solver on a large mesh composed of
over 6 million tetrahedrons (precisely \num{6533341}).
The purpose of this test is to validate
the scaling of the solver,
as well as display its results on a realistic test case.
The treated mesh represents an anthropomorphic mannequin
named Kyoto, and it has been used in other works.
So far, it has not been included in peer-reviewed articles,
but it was used in several PhD theses
\cite{weber:tel-01911261,Hou2020}
and in a PRACE SHAPE project in collaboration with
the AxesSim company,
see the white paper \cite{GirWebCirCam2018}
and the report \cite{GENCI2020}.
In these contexts, the tetrahedra were cut into  four
hexahedra each,
which increased the total number of elements fourfold.

The body model is composed of 12 organs
(including the skeleton and the skin).
To handle these different body parts,
and to include a source term modeling
the behavior of a current generated by an antenna,
we modify Maxwell's equations \eqref{eq:Maxwell_in_PDE_form},
as follows:
\begin{equation}
	\label{eq:Maxwell_with_eps_in_PDE_form}
	\begin{dcases}
		\partial_t (\varepsilon E) - \nabla \times H = - \sigma E - J, \\
		\partial_t H + \nabla \times E = 0.
	\end{dcases}
\end{equation}
In \eqref{eq:Maxwell_with_eps_in_PDE_form},
$J$ is the time- and space-dependent electric current density and
$\varepsilon$ is the permittivity of the material,
which obviously depends on the material
(and therefore on the space variable $X$).
We introduce the vacuum permittivity $\varepsilon_0$,
to write $\varepsilon = \varepsilon_r \varepsilon_0$
with $\varepsilon_r$
the relative permittivity of the material.
For each body part, the values of $\varepsilon_r$ and $\sigma$
are listed in \cref{tab:materials_kyoto}.
In addition, we set $\widetilde E = \varepsilon_r E$.
Reformulating \eqref{eq:Maxwell_with_eps_in_PDE_form},
we obtain, assuming that each material has the same permeability:
\begin{equation}
	\label{eq:Maxwell_with_eps_in_flux_in_PDE_form}
	\begin{dcases}
		\partial_t \widetilde E - \nabla \times H = - \sigma \frac{\widetilde E}{\varepsilon_r} - \frac{1}{\varepsilon_r} J, \\
		\partial_t H + \nabla \times \frac{\widetilde E}{\varepsilon_r} = 0.
	\end{dcases}
\end{equation}
Note that the flux in
\eqref{eq:Maxwell_with_eps_in_flux_in_PDE_form}
is discontinuous as soon as $\varepsilon_r$ is discontinuous,
which is the case here since different materials
have different relative permittivities.
In practice, we solve for $\widetilde E$ and $H$.

\begin{table}[!ht]
	\centering
	\caption{%
		Electromagnetic constants for each material
		in the anthropomorphic mannequin simulation.
	}
	\label{tab:materials_kyoto}
	\begin{tabular}{lcc}
		\toprule
		Material    & $\varepsilon_r$ & $\sigma$ (Sv \, m\textsuperscript{-1}) \\ \cmidrule(lr){1-3}
		brain       & 48.34           & 2.02                                 \\
		heart       & 58.67           & 3.02                                 \\
		lung        & 22              & 0.36                                 \\
		liver       & 41.82           & 1.9                                  \\
		gallbladder & 60              & 2                                    \\
		spleen      & 56.75           & 2.46                                 \\
		pancreas    & 56.75           & 2.46                                 \\
		kidney      & 56.83           & 2.62                                 \\
		colon       & 48.5            & 0.93                                 \\
		bladder     & 20              & 0.7                                  \\
		muscle      & 50              & 1.33                                 \\
		bone        & 11.41           & 0.43                                 \\
		cartilage   & 36              & 1.6                                  \\ \bottomrule
	\end{tabular}
\end{table}

A volume-meshed dipole antenna has been placed
next to the left arm of the body model,
and a volumetric source term has been imposed along that dipole.
This leads to the electric current density
$J$ in \eqref{eq:Maxwell_with_eps_in_flux_in_PDE_form}
being nonzero only within the antenna.
In dimensional quantities,
the cells in the antenna have average size
$2 \cdot 10^{-4}$ m,
while the cells in vacuum have average size
$2 \cdot 10^{-2}$ m.
This means that the ratio of largest cell size
over smallest cell size is around $100$;
hence, this mesh provides a good framework
to test our CFL-less methodology.
In addition, to mimic a Bluetooth antenna,
this dipole antenna emits a
modulated Gaussian pulse,
with dimensional frequency 2.4 GHz,
which lasts for about 1.5 ns.
After that time has elapsed, the source term $J$
vanishes in the whole domain.

\begin{table}[!ht]
    \centering
	\caption{Differences between CLAC and KOUGLOFV.}
	\label{tab:differences_CLAC_KOUGLOFV}
	\begin{tabular}{rcc}
		\toprule
		                         & CLAC                                     & KOUGLOFV                                 \\ \cmidrule(lr){1-3}
		space order              & 3 (10-point tetrahedra)                  & 3 (10-point tetrahedra)                  \\
		time order               & 3 (explicit, RK3)                        & 2 (implicit, Crank-Nicolson)             \\
		floating-point precision & simple                                   & double                                   \\
		boundary conditions      & Silver-Müller                            & homogeneous Dirichlet                    \\
		time step                & $\Delta t = 9.04 \cdot 10^{-6}$ ns & $\Delta t = 3.33 \cdot 10^{-4}$ ns \\
		\bottomrule
	\end{tabular}
\end{table}

In \cref{tab:differences_CLAC_KOUGLOFV},
we sum up the differences between the two solvers.
Since the two codes have
different time stepping strategies,
different boundary conditions, etc.,
we do not expect their behaviors
to be quantitatively comparable.
However, they should, qualitatively,
lead to similar results.
These differences between the two solvers also explain why
the KOUGLOFV code takes about as much computation time on
a few hundred CPUs than the CLAC code on 6 GPUs.
The computation time figures are reported in
\cref{tab:CPU_times_kyoto},
and we observe that the KOUGLOFV code takes about as much time
to run as the CLAC solver despite being in double-precision
arithmetic and running on $128$ CPU cores rather than $6$ GPUs.

\begin{table}
    \centering
	\caption{%
		Computation time for both KOUGLOFV and CLAC codes,
		for the anthropomorphic mannequin simulation.
	}
	\label{tab:CPU_times_kyoto}
	\begin{tabular}{ccc}
		\toprule
		solver   & hardware platform & computation time  \\
		\cmidrule(lr){1-3}
		CLAC     & \makecell{
		6 $\times$ Nvidia GeForce GTX 1080 Ti            \\
			single precision arithmetic
		}        & 31 h                       \\[\bigskipamount]
		KOUGLOFV & \makecell{
		AMD EPYC 7713x2, 128 cores, 2 GHz \\
			double precision arithmetic
		}        & 20 h                       \\
		\bottomrule
	\end{tabular}
\end{table}

Indeed, this good behavior happens thanks to the CFL-less
implicit time-stepping in the scheme underlying in the KOUGLOFV solver.
To correctly implement this implicit scheme, we have to give a value to
the over-relaxation parameter $\omega \in [1, 2]$ from \eqref{eq:over_relaxation},
see the discussion at the end of \cref{sec:kinetic_algorithm}.
Recall that taking $\omega = 2$ leads to a second-order scheme,
while other choices lead to a first-order scheme,
with larger values of $\omega$ corresponding to higher resolutions.
Here, we made the choice to take $\omega = 1.8$,
since larger values led to small instabilities stemming
from the discontinuous flux function in \eqref{eq:Maxwell_with_eps_in_flux_in_PDE_form}:
indeed, lowering the order of the scheme helped curb these spurious oscillations.
This is reminiscent of ideas from, for instance, \cite{HigHapKocKup2014,MicTho2022}.
In addition, another choice to make is the value of the time step,
since the scheme is unconditionally stable.
We chose to take $\Delta t = \num{1e-4}$ in non-dimensional form,
leading to a time step about $40$ times larger than the stability limit of
the scheme from the CLAC code.
This made it possible to have good results from a relatively short simulation.

\begin{figure}[!ht]
	\centering
	\includegraphics[width=0.45\textwidth]{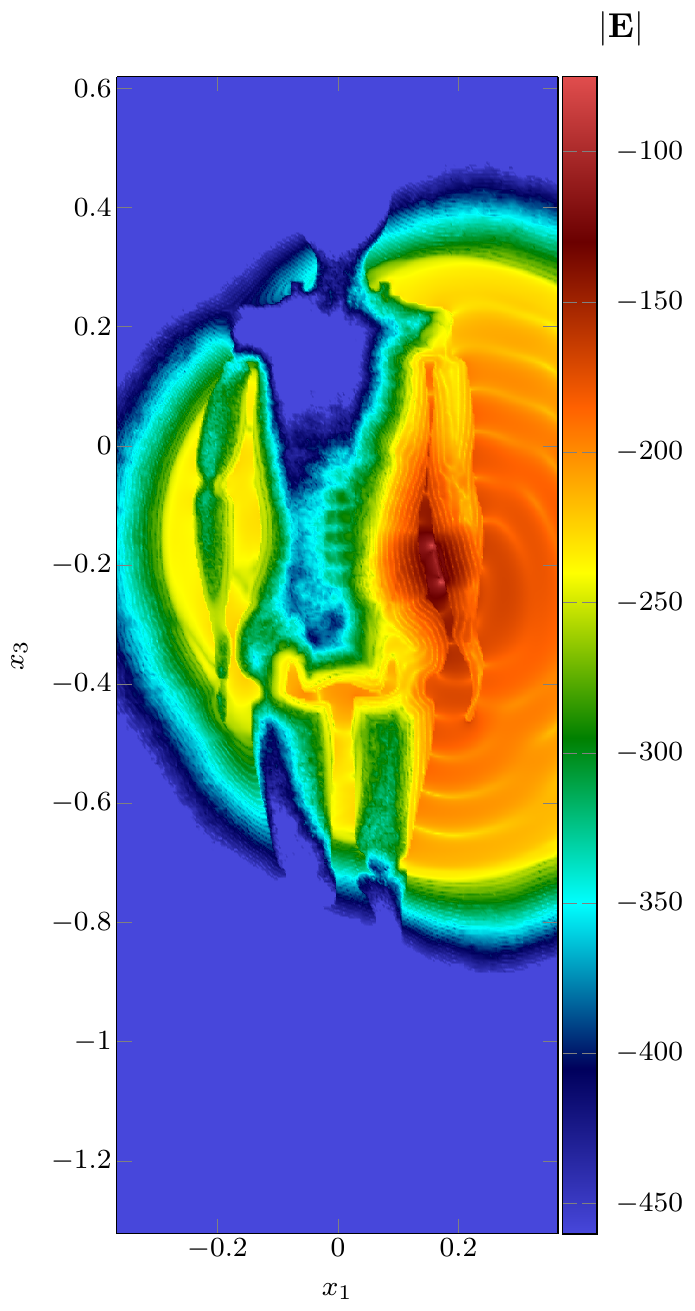}
	\includegraphics[width=0.45\textwidth]{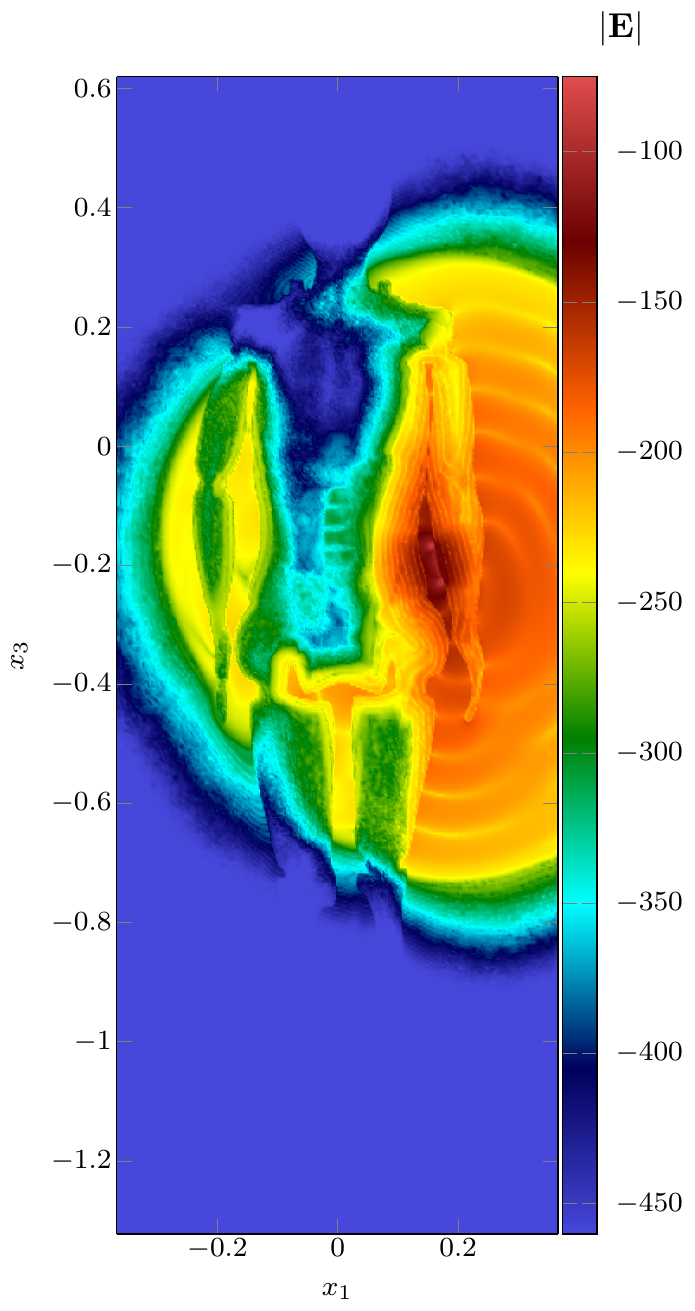}
	\caption{%
		Numerical approximation of $| E |$,
        sliced in the $(x_1, x_3)$ plane,
		at $t = \SI{0.5}{\nano\second}$.
		Left panel: results from CLAC;
		right panel: results from KOUGLOFV.
	}
	\label{fig:kyoto_t_0v5}
\end{figure}

\begin{figure}[!ht]
	\centering
	\includegraphics[width=0.45\textwidth]{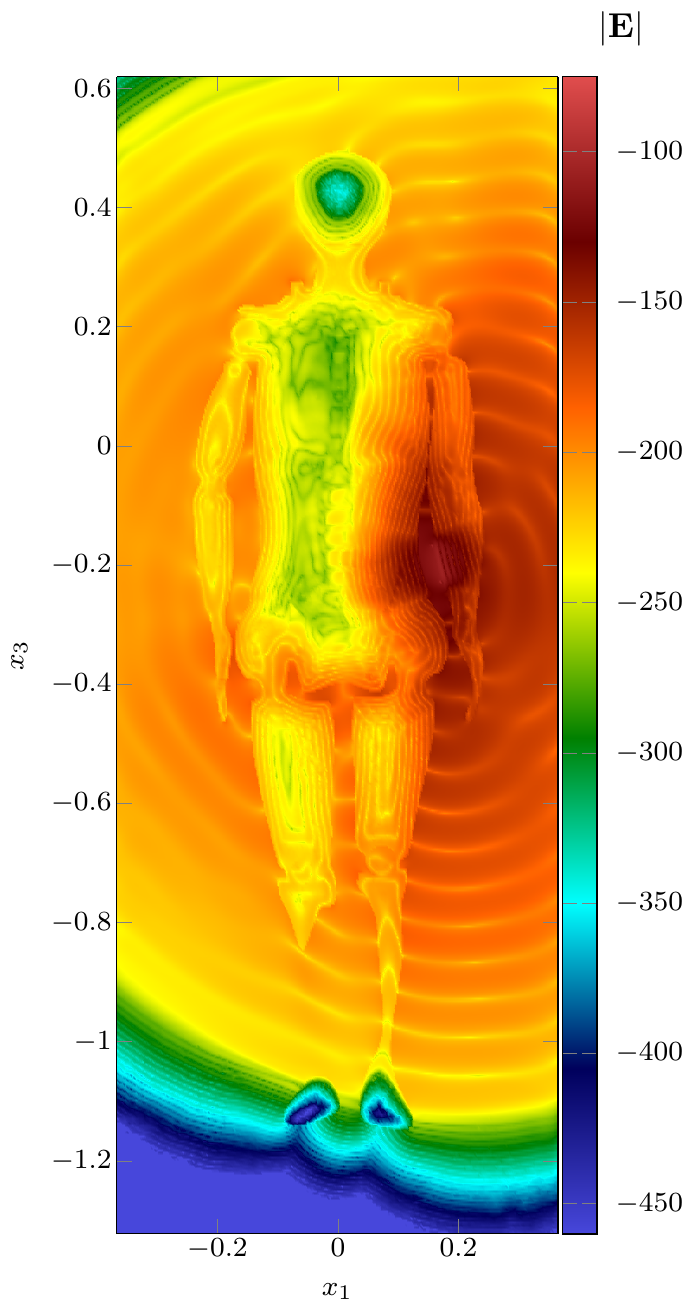}
	\includegraphics[width=0.45\textwidth]{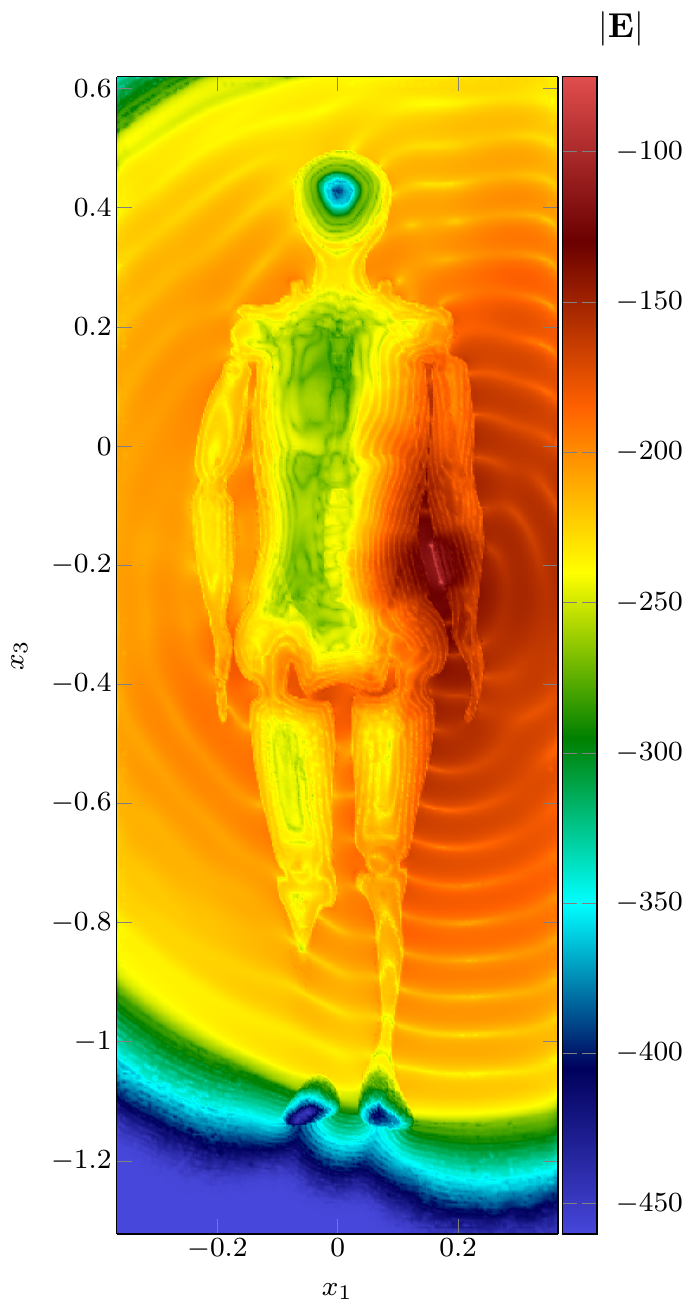}
	\caption{%
		Numerical approximation of $| E |$,
        sliced in the $(x_1, x_3)$ plane,
		at $t = \SI{0.9}{\nano\second}$.
		Left panel: results from CLAC;
		right panel: results from KOUGLOFV.
	}
	\label{fig:kyoto_t_0v9}
\end{figure}

\begin{figure}[!ht]
	\centering
	\includegraphics[width=0.45\textwidth]{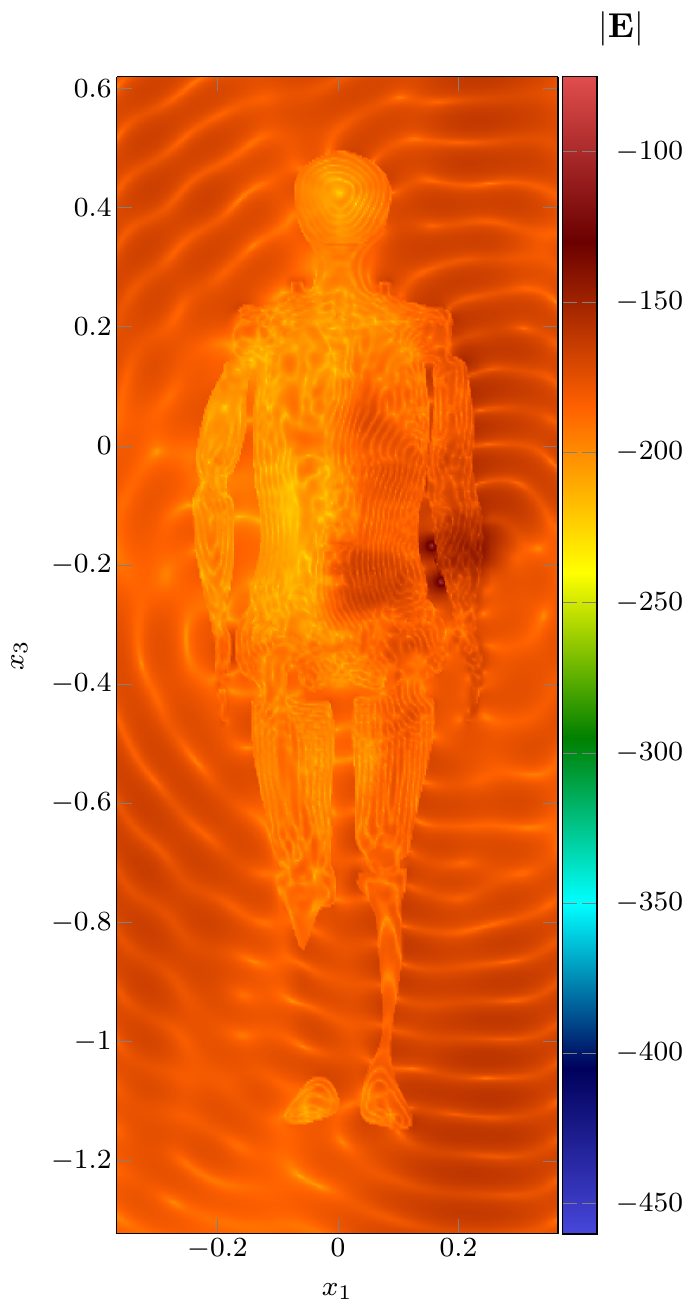}
	\includegraphics[width=0.45\textwidth]{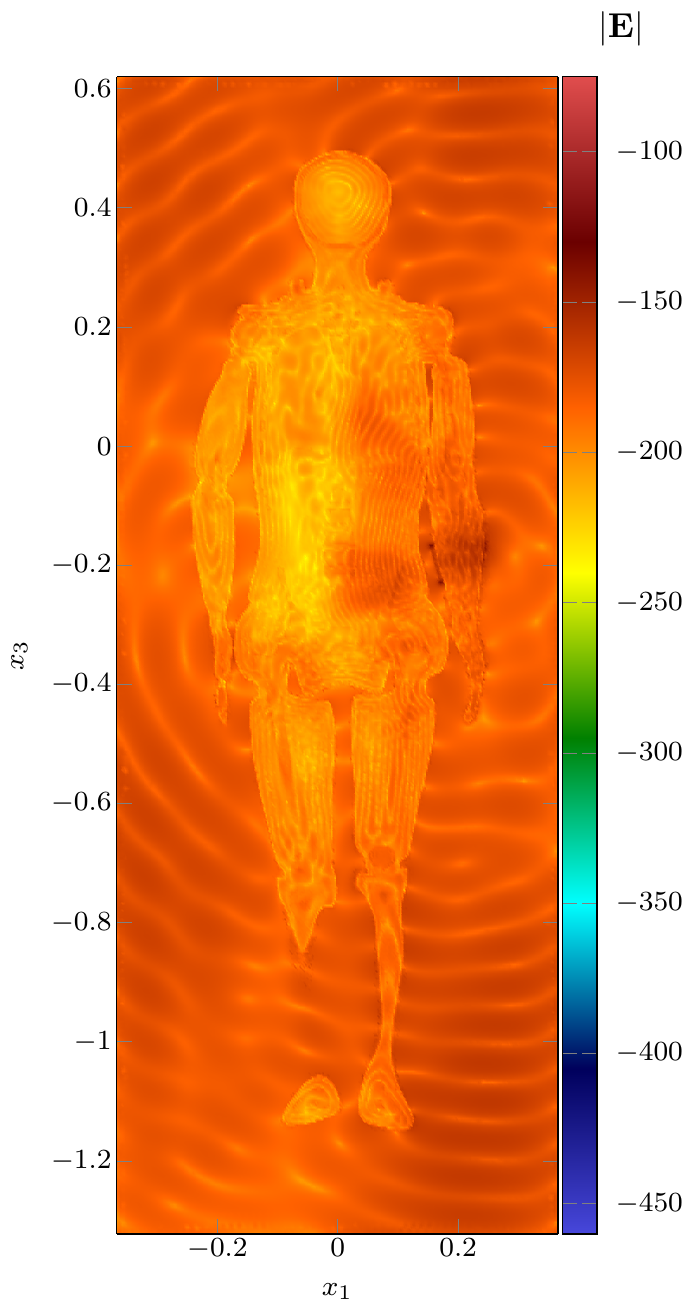}
	\caption{%
		Numerical approximation of $| E |$,
        sliced in the $(x_1, x_3)$ plane,
		at $t = \SI{1.8}{\nano\second}$.
		Left panel: results from CLAC;
		right panel: results from KOUGLOFV.
	}
	\label{fig:kyoto_t_1v8}
\end{figure}

We first present the results in the $(x_1, x_3)$ plane.
They are displayed on \cref{fig:kyoto_t_0v5}
for $t = 0.5$ ns,
on \cref{fig:kyoto_t_0v9}
for $t = 0.9$ ns,
and on \cref{fig:kyoto_t_1v8}
for $t = 1.8$ ns.
In each case, we observe good agreement between
the results of CLAC and KOUGLOFV,
despite the differences in the two approaches.

\begin{figure}[!ht]
	\centering
	\includegraphics{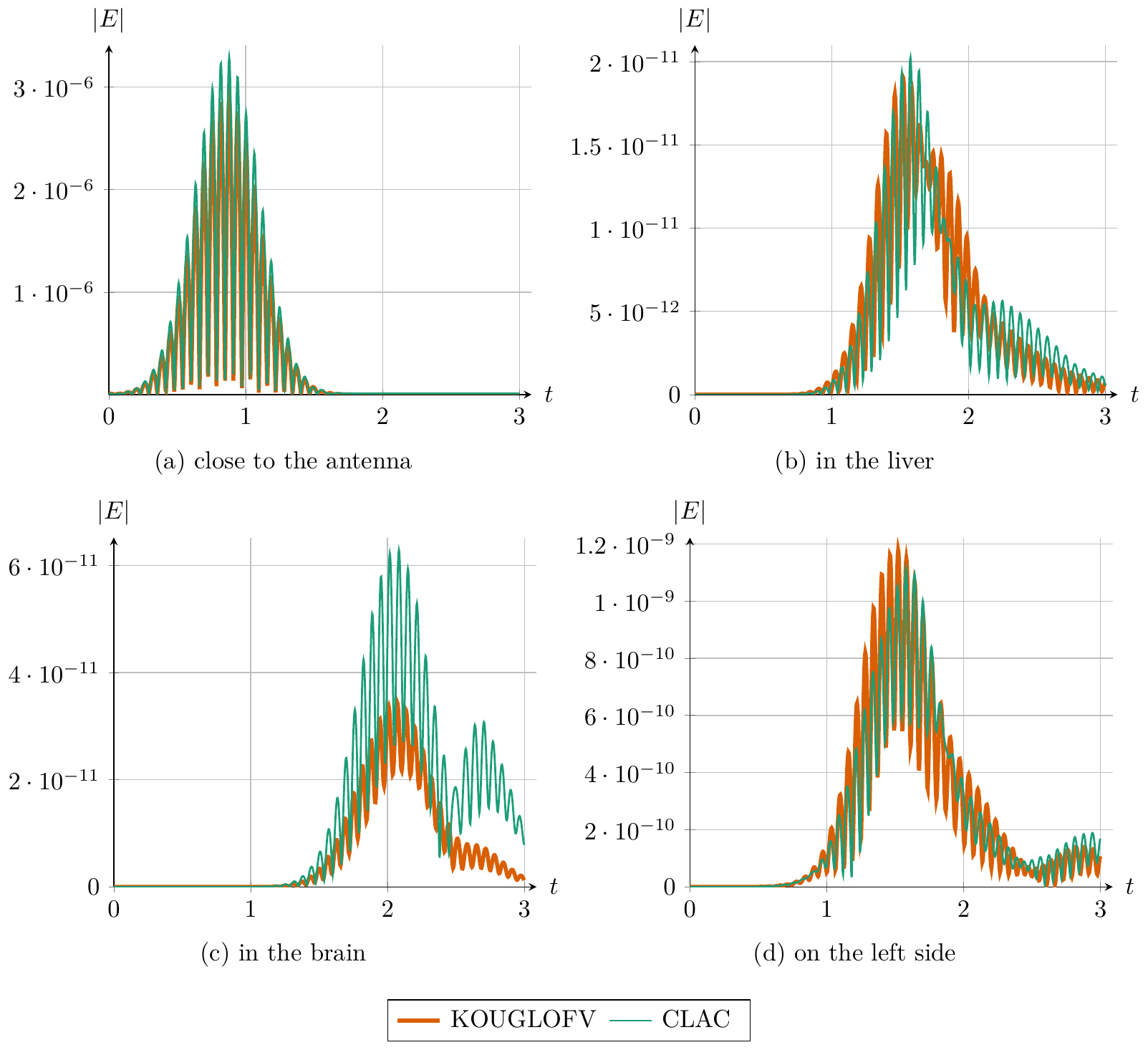}
	\caption{%
		Numerical solution at four points within the mesh:
		comparison between the solutions given by the two codes.
	}
	\label{fig:kyoto_comparison_at_points}
\end{figure}

\begin{figure}[!ht]
	\centering
	\includegraphics{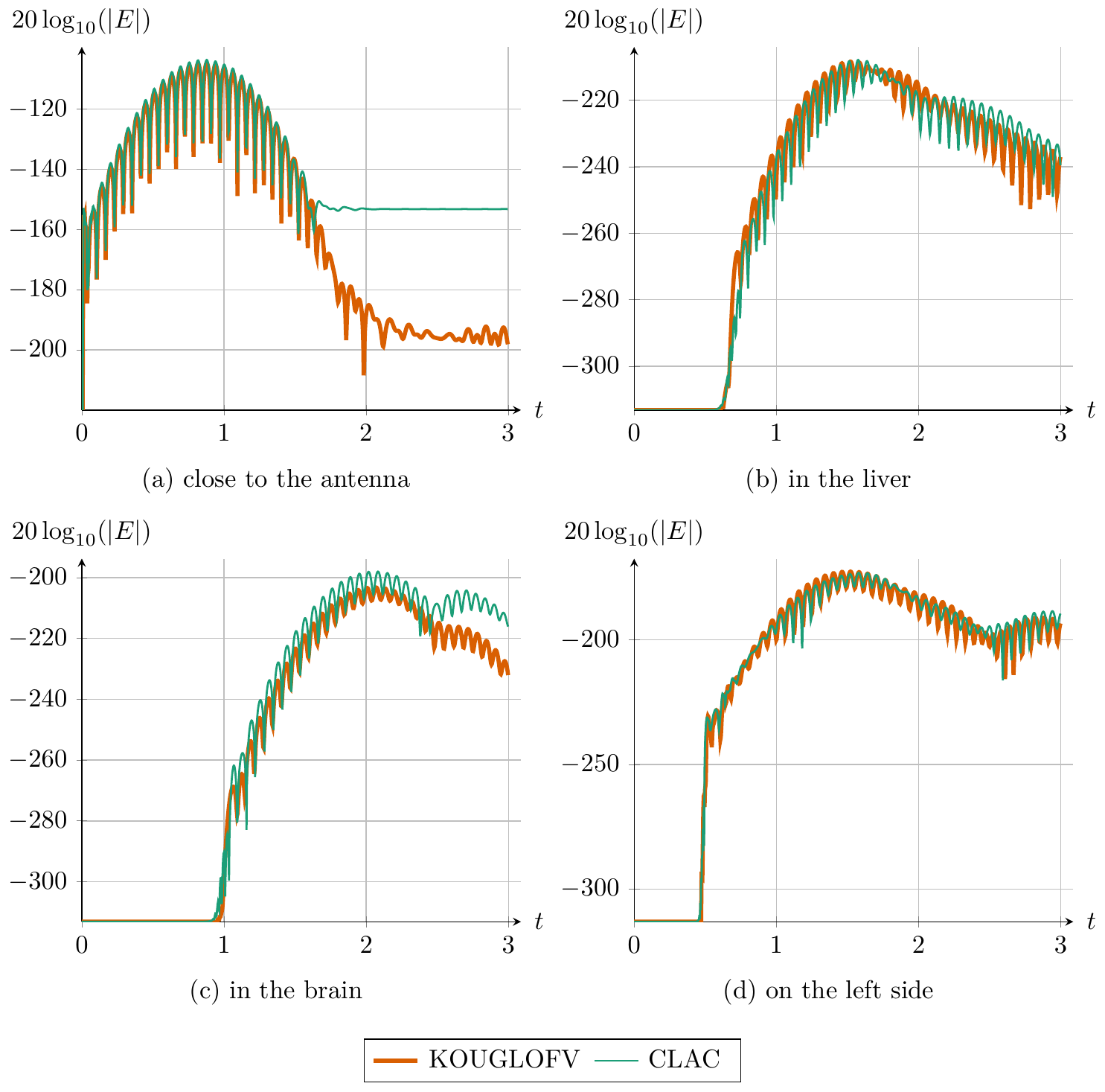}
	\caption{%
		Numerical solution at four points within the mesh:
		comparison between the solutions given by the two codes.
        Here, the results are presented in dB,
        corresponding to taking $20$ times the base-$10$
        logarithm of $|E|$.
	}
	\label{fig:kyoto_comparison_at_points_dB}
\end{figure}

\begin{figure}[!ht]
	\centering
	\includegraphics{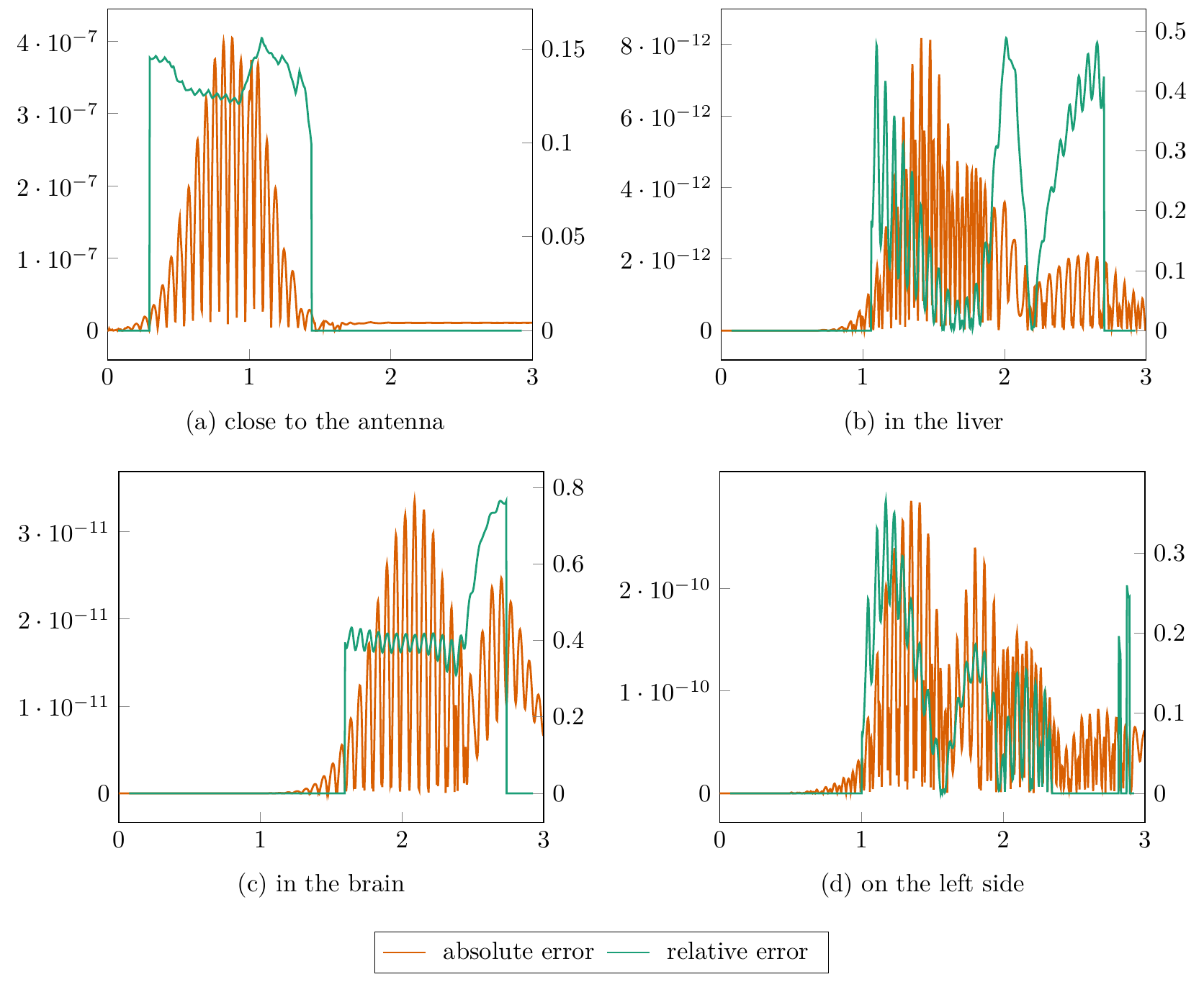}
	\caption{%
		Absolute (left axes) and relative (right axes) errors
		between the KOUGLOFV and CLAC solvers,
		at four points within the mesh.
	}
	\label{fig:kyoto_comparison_errors}
\end{figure}

To get more precise results, we compare in
\cref{fig:kyoto_comparison_at_points,fig:kyoto_comparison_at_points_dB,fig:kyoto_comparison_errors}
the two numerical results at four points:
\begin{enumerate}[nosep]
	\item close to the antenna (in vacuum), top left panels;
	\item in the liver, top right panels;
	\item in the brain, bottom left panels;
	\item on the left side (in vacuum), bottom right panels.
\end{enumerate}
At each point, the signals have the same shape for the two codes,
although the signal from KOUGLOFV is more diffused
compared to the one from CLAC.
Close to the antenna, the signal produced by KOUGLOFV
has about 15\% relative error with respect to the CLAC signal,
which can be attributed to the difference in source term
discretizations between the two implementations.
That relative error can be used as a baseline to compare the
two results.
Moving further away from the source, for instance in the brain,
the two signals have the same shape but the signal
from KOUGLOFV is more diffused than the one from CLAC:
this can be attributed to the fact that lowering the order of
the scheme was necessary to handle the discontinuous flux
in this simulation.

\section{Conclusion}
\label{sec:conclusion}

We presented an adaptation of the kinetic DG method introduced in
\cite{gerhard2022unconditionally}. The method can handle arbitrary
conservation laws and complex unstructured meshes. It has the complexity
of a time-explicit scheme but is CFL-free.

The method presents good parallelization features, for both shared
memory and distributed memory computers. To improve the parallel
scaling on distributed memory computers, we have proposed a subdomain
decomposition method that relaxes the task dependencies of the kinetic
scheme but keeps the possibility to use large time steps. The method
has been tested and validated on realistic electromagnetic simulations.

In our future works, we plan to apply the method to other conservation
laws arising for instance in the modeling of multiphase compressible
flows. Investigations are also needed for a more rigorous treatment
of the boundary conditions.

\bibliographystyle{plain}
\bibliography{kinetic-cfl-less}

\end{document}